\documentclass[12pt,a4paper]{article}
\usepackage{latexsym,epsfig,fullpage,amsfonts,amssymb,amsmath,amscd,epic,graphics,rotating,theorem}
\usepackage[all]{xy}
\usepackage[usenames]{color}
\baselineskip 18pt \textwidth 15cm  \theoremstyle{plain}
\newtheorem{lemma}{Lemma}[section]

\newtheorem{remark}[lemma]{Remark}

\newtheorem{theorem}{Theorem}[section]

\newtheorem{corollary}{Corollary}

{\theorembodyfont{\rmfamily}
\begin{document}
\baselineskip=18pt
\newcommand{\Appendix}[2][?]{%
\refstepcounter{section}%
%\addcontentsline{toc}{appendix}%
{\protect\numberline{\appendixname~\thesection}#1}%
{\flushright\large\bfseries\appendixname\ \thesection\par
\nonhyphens\centering#2\par}%
\sectionmark{#1}\vspace{\baselineskip}}
\newcommand{\sAppendix}[1]{%
{\flushright\large\bfseries\appendixname\par
\nonhyphens\centering#1\par}%
\vspace{\baselineskip}}
\newcommand{\N}{{\mathbb N}}\newcommand{\Q}{{\mathbb Q}}
\newcommand{\oz}{{\overline z}}
\newcommand{\PP}{{\mathbb P}}
\newcommand{\Z}{{\mathbb Z}}\newcommand{\A}{{\mathbb A}}
\newcommand{\R}{{\mathbb R}}
\newcommand{\C}{{\mathbb C}}
\newcommand{\K}{{\mathbb K}}
\newcommand{\F}{{\mathbb F}}\newcommand{\bbS}{{\mathbb S}}
\newcommand{\proofend}{\hfill$\Box$\bigskip}
\newcommand{\eps}{{\varepsilon}}
\newcommand{\ko}{{\mathcal O}}
\newcommand{\ox}{{\overline{\boldsymbol x}}}
\newcommand{\wz}{{\widetilde z}}
\newcommand{\wa}{{\widetilde a}}
\newcommand{\bz}{{\boldsymbol z}}
\newcommand{\bp}{{\boldsymbol p}}\newcommand{\by}{{\boldsymbol y}}
\newcommand{\wy}{{\widetilde y}}
\newcommand{\wc}{{\widetilde c}}
\newcommand{\bi}{{\omega}}
\newcommand{\bx}{{\boldsymbol x}}\newcommand{\bv}{{\boldsymbol v}}\newcommand{\bw}{{\boldsymbol w}}
\newcommand{\Log}{{\operatorname{Log}}}
\newcommand{\Graph}{{\operatorname{Graph}}}
\newcommand{\jet}{{\operatorname{jet}}}\newcommand{\red}{{\operatorname{red}}}
\newcommand{\Tor}{{\operatorname{Tor}}}\newcommand{\defn}{d}
\newcommand{\sqh}{{\operatorname{sqh}}}
\newcommand{\const}{{\operatorname{const}}}
\newcommand{\Arc}{{\operatorname{Arc}}}\newcommand{\ini}{{\operatorname{Ini}}}
\newcommand{\Sing}{{\operatorname{Sing}}}
\newcommand{\Span}{{\operatorname{Span}}}
\newcommand{\Aut}{{\operatorname{Aut}}}
\newcommand{\Ker}{{\operatorname{Ker}}}
\newcommand{\Int}{{\operatorname{Int}}}
\newcommand{\Aff}{{\operatorname{Aff}}}
\newcommand{\Area}{{\operatorname{Ar}}}
\newcommand{\val}{{\operatorname{Val}}}
\newcommand{\conv}{{\operatorname{Conv}}}
\newcommand{\conj}{{\operatorname{Conj}}}\newcommand{\grad}{{\operatorname{grad}}}
\newcommand{\Iso}{{\operatorname{Iso}}}
\newcommand{\rk}{{\operatorname{rk}}}
\newcommand{\pr}{{\operatorname{pr}}}
\newcommand{\op}{{\overline{\boldsymbol p}}}
\newcommand{\ov}{{\overline v}}
\newcommand{\ks}{{\cal S}}
\newcommand{\kc}{{\cal C}}
\newcommand{\ki}{{\cal I}}
\newcommand{\kj}{{\cal J}}
\newcommand{\ke}{{\cal E}}
\newcommand{\kz}{{\cal Z}}
\newcommand{\tet}{{\theta}}
\newcommand{\Del}{{\Delta}}
\newcommand{\bet}{{\beta}}
\newcommand{\mm}{{\mathfrak m}}
\newcommand{\kap}{{\kappa}}
\newcommand{\del}{{\delta}}
\newcommand{\sig}{{\sigma}}
\newcommand{\alp}{{\alpha}}
\newcommand{\Sig}{{\Sigma}}
\newcommand{\Gam}{{\Gamma}}
\newcommand{\gam}{{\gamma}}
\newcommand{\Lam}{{\Lambda}}
\newcommand{\lam}{{\lambda}}
\title{Welschinger invariants of toric Del Pezzo surfaces with non-standard real structures}
\author{E. Shustin}
%\thanks{The
%author has been supported by the Hermann-Minkowski Minerva Center
%for Geometry at the Tel Aviv University.}}
\date{}
\maketitle
\begin{abstract} The Welschinger invariants of real rational algebraic surfaces
are natural analogues of the Gromov-Witten invariants, and they
estimate from below the number of real rational curves passing
through prescribed configurations of points. We establish a
tropical formula for the Welschinger invariants of four toric Del
Pezzo surfaces, equipped with a non-standard real structure. Such
a formula for real toric Del Pezzo surfaces with a standard real
structure (i.e., naturally compatible with the toric structure)
was established by Mikhalkin and the author. As a consequence we
prove that, for any real ample divisor $D$ on a surfaces $\Sig$
under consideration, through any generic configuration of
$c_1(\Sig)D-1$ generic real points there passes a real rational
curve belonging to the linear system $|D|$.
\end{abstract}

\section*{Introduction}

The Welschinger invariants \cite{Wel,Wel1} play a central role in
the enumerative geometry of real rational curves on real rational
surfaces, providing lower bounds for the number of real rational
curves passing through generic, conjugation invariant
configurations of points, whereas the number of respective complex
curves (Gromov-Witten invariant) serves as an upper bound. Methods
of the tropical enumerative geometry, developed in
\cite{M2,M3,ShP}, allowed one to compute and estimate the
Welschinger invariants for the real toric Del Pezzo surfaces,
equipped with the standard real structure \cite{IKS,IKS2,ShW}: the
plane $\PP^2$, the plane $\PP^2_k$ with blown up $k=1$, $2$, or
$3$ real (generic) points, and the quadric $(\PP^1)^2$.

We notice that the available technique of the tropical enumerative
geometry applies only to toric surfaces, and among them the
Welschinger invariant is well-defined only for
unndodal\footnote{Like in \cite{IKS,IKS2} "unnodal" means the
absence of $(-n)$-curves, $n\ge 2$.} Del Pezzo surfaces. So, the
main goal of this paper is to compute Welschinger invariants for
other real toric Del Pezzo surfaces using the tropical enumerative
geometry.

Along the Comessatti's classification of real rational surfaces
\cite{C1,C2} (see also \cite{K1}), besides the standard real toric
Del Pezzo surfaces, there are five more types, having a non-empty
real points set, which we call {\it non-standard} and denote as
$\bbS^2$, the quadric whose real point set is a sphere,
$\bbS^2_{1,0},\bbS^2_{2,0},\bbS^2_{0,2}$, the sphere with blown up
one or two real points, or a pair of conjugate imaginary points,
respectively, and, at last, $(\PP^1)^2_{0,2}$, the standard real
quadric blown up at two imaginary conjugate points.

In the present paper we derive the tropical formula for the
Welschinger invariants of $\bbS^2$, $\bbS^2_{1,0}$,
$\bbS^2_{2,0}$, and $\bbS^2_{0,2}$\footnote{On the real toric Del
Pezzo surfaces with empty real point set there are no real
rational curves, and thus, Welschinger invariants vanish.}, that
is we express the Welschinger invariants as the sums of weights of
certain discrete combinatorial objects, running over specified
finite sets, and which are related to tropical curves
corresponding to the real algebraic curves in count. The surface
$(\PP^1)^2_{0,2}$ requires a completely different treatment, and
it will be addressed in a forthcoming paper.

The formulation of our result is split into Theorem \ref{tsp1},
section \ref{secsp3} (all the surfaces and totally real
configurations points), Theorem \ref{tsp3}, section \ref{secspn1}
(surfaces $\bbS^2$, $\bbS^2_{1,0}$, $\bbS^2_{2,0}$ and
configurations with imaginary points), and Theorem \ref{tsp4},
section \ref{secspn4} (surface $\bbS^2_{0,2}$ and configurations
with imaginary points). The reason to have a separate statement
for the case of totally real configurations is a simpler
formulation and its special importance in applications. For
example, we prove the positivity of the Welschinger invariants,
related to the totally real configurations, which immediately
implies the existence of real rational curves belonging to given
linear systems and passing through generic configurations of
suitable number of real points (see Corollary \ref{csp2}, section
\ref{secsp3}).

As compared with the standard real toric Del Pezzo surfaces
\cite{M3,ShP,ShW}, the case considered in the present paper looks
much more degenerate, for instance, the plane tropical curves
(amoebas) in count are highly reducible. On the other hand,
similarly to the standard case, the final answer is basically
expressed as the weighted number of irreducible rational
parameterizations of the above plane tropical curves, though the
parameterizations generically are not trivalent.

Finally, we notice that, in a joint work with I. Itenberg and V.
Kharlamov \cite{IKS4}, using Theorem \ref{tsp1} we prove the
asymptotic relation $$\lim_{n\to\infty}\frac{\log W_0(\Sig,{\cal
L}^{\otimes n})}{n\log n}=\lim_{n\to\infty}\frac{\log
N_0(\Sig,{\cal L}^{\otimes n})}{n\log n}=-c_1({\cal L})K_\Sig\ ,$$
for any real ample line bundle ${\cal L}$ on a surface
$\Sig=\bbS^2$, $\bbS^2_{1,0}$, $\bbS^2_{2,0}$, or $\bbS^2_{0,2}$,
which before was established for the standard real toric Del Pezzo
surfaces \cite{IKS2}.

\medskip
\noindent {\bf Welschinger invariants.} For the reader's
convenience, we recall the definition of Welschinger invariants.
Let $\Sig$ be a real toric Del Pezzo surface with a non-empty real
part, ${\cal L}$ a very ample real line bundle on $\Sig$, and let
non-negative integers $r',r''$ satisfy
\begin{equation}r'+2r''=-c_1({\cal L})K_\Sig-1\
.\label{enn36}\end{equation} Denote by $\Omega_{r''}(\Sig,{\cal
L})$ the set of configurations of $-c_1({\cal L})K_\Sig- 1$
distinct points of $\Sig$ such that $r'$ of them are real and the
rest consists of $r''$ pairs of imaginary conjugate points. The
Welschinger number $W_{r''}(\Sig,{\cal L})$ is the sum of weights
of all the real rational curves in the linear system $|{\cal L}|$,
passing through a generic configuration
$\op\in\Omega_{r''}(\Sig,{\cal L})$, where the weight of a real
rational curve $C$ is $1$ if it has an even number of real
solitary nodes, and is $-1$ otherwise. Since the complex structure
of $\Sig$ determines a symplectic structure, which is generic in
its deformation class, by Welschinger's theorem \cite{Wel,Wel1},
$W_{r''}(\Sig,{\cal L})$ does not depend on the choice of a
generic element $\op\in\Omega_{r''}(\Sig,{\cal L})$ (a simple
proof of this fact can be found in \cite{IKS1}). The definition
immediately implies the inequality
\begin{equation}|W_{r''}(\Sig,{\cal L})|\le R_{\Sig,{\cal L}}(\op)\le N_{\Sig,{\cal L}}\
,\label{enn200}\end{equation} where $R_{\Sig,{\cal L}}(\op)$ is
the number of real rational curves in $|{\cal L}|$ passing through
a generic configuration $\op\in\Omega_{r''}(\Sig,{\cal L})$, and
$N_{\Sig,{\cal L}}$ is the number of complex rational curves in
$|{\cal L}|$, passing through generic $-c_1({\cal L})K_\Sig- 1$
points in $\Sig$.

\medskip
\noindent {\bf A tropical calculation of the Welschinger
invariant.} Our approach to calculating the Welschinger invariants
is quite similar to that in \cite{IKS,ShW}, and it heavily relies
on the enumerative tropical algebraic geometry developed in
\cite{M2,M3,ShP}. More precisely, we replace the complex field
$\C$ by the field $\K=\bigcup_{m\ge 1}\C\{\{t^{1/m}\}\}$ of the
complex, locally convergent Puiseux series endowed with the
standard complex conjugation and with a non-Archimedean valuation
$$\val:\K^*\to\R,\quad\val\left(\sum_ka_kt^k\right)=-\min\{k\ :\ a_k\ne 0\}\
.$$ A rational curve over $\K_\R$, belonging to a linear system
$|{\cal L}|_\K$ and passing through a generic configuration
$\op\in\Omega_{r''}(\Sig(\K),{\cal L})$, is viewed as an
equisingular family of real rational curves in $\Sig$ over the
punctured disc. We construct an appropriate limit of the family of
surfaces and embedded curves at the disc center. The central
surface is usually reducible, and the adjacency of its components
is encoded by a tropical curve in the real plane, which passes
through the configuration $\val(\op)\subset\R^2$. The central
curve is split into components called {\it limit curves}. The pair
(tropical curve, limit curves) is called the {]it tropical limit}
of the given curve $C\in|{\cal L}|_\K$.

We precisely describe the tropical limits of real rational curves
passing through generic configurations of real points in
$\Sig(\K)$, then compute the {\it Welschinger weights} of the
respective tropical curves, i.e., the contribution to the
Welschinger invariant of the real algebraic curves projecting to
the given tropical curve. The result, accumulated in Theorem
\ref{tsp1} (section \ref{secsp1}), represents the Welschinger
invariants $W_0(\Sig,{\cal L})$ as the numbers of some
combinatorial objects, forming finite discrete sets.

The proof is based on the techniques of \cite{ShP,ShW}, both in
the determining tropical limits and in the patchworking
construction, which recovers algebraic curves over $\K$ from their
tropical limits. We should like to remark that the answer rather
differs from that for the standard real Del Pezzo surfaces.
Namely, in the standard case, the tropical limits are basically
encoded by tropical curves, which are rational and irreducible. In
the non-standard case, one obtains a relatively small number of
possible tropical curves, which split into unions of some
primitive tropical curves. In contrast, the weights of the
tropical curves are large and are defined in a non-trivial
combinatorial way.

We also notice that the patchworking theorems from \cite{ShP,ShW}
cover our needs in the present paper. In contrast, the
determination of tropical limits meets extra difficulties, caused
by the fact that the {\it generic} configurations of real points
on the surfaces under consideration project by
$\val:(\K^*)^2\to\R^2$ to {\it non-generic} configurations in
$\R^2$ (cf. a similar problem in \cite{ShP}).

\medskip
\noindent{\bf Applications to enumerative geometry.} From Theorem
\ref{tsp1} we immediately derive the positivity of the Welschinger
invariants in the considered situations, which in view of
(\ref{enn200}) results in Corollary \ref{csp2}, section
\ref{secsp1}, which says that , for any real very ample line
bundle ${\cal L}$ on a non-standard real toric Del Pezzo surface
$\Sig=\bbS^2$, $\bbS^2_{1,0}$, $\bbS^2_{2,0}$, or $\bbS^2_{0,2}$
and any generic configuration of $-c_1({\cal L})K_\Sig-1$ real
points on $\Sig$ there exists a real rational curve $C\in|{\cal
L}|$ passing through the given configuration.

\medskip\noindent
{\bf Acknowledgements.} The research was supported by the grant
no. 465/04 from the Israel Science Foundation, by the grant from
the High Council for Scientific and Technological Cooperation
between France and Israel, and by the Hermann-Minkowski Minerva
center for Geometry at the Tel Aviv University. A part of this
work was done during the author's stay at the Max-Planck-Institut
f\"ur Mathematik (Bonn) and at the Institut Henri Poincar\'e
(Paris). The author is very grateful to MPI and IHP for the
hospitality and excellent work conditions. I would like to thank
I. Itenberg and V. Kharlamov for useful discussions, which were
vital for completing this work.

\section{Tropical formula for the Welschinger
invariants}\label{secspn9}

\subsection{Lattice polygons associated with the non-standard real toric Del Pezzo surfaces}
\label{sec1} 
\begin{figure}
\begin{center}
\epsfxsize 145mm \epsfbox{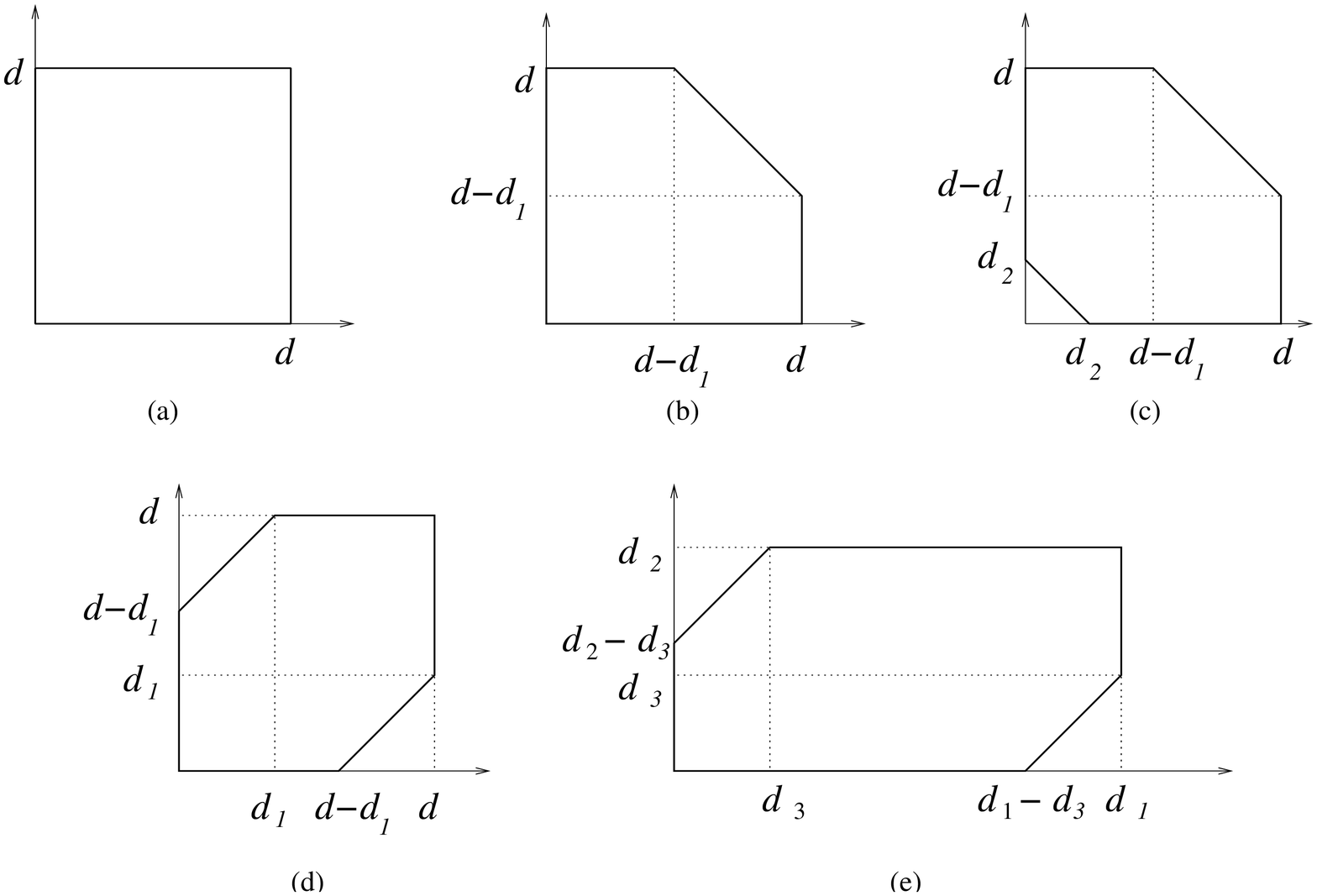}
\end{center}
\caption{Polygons associated with $\bbS^2$, $\bbS^2_{1,0}$,
$\bbS^2_{2,0}$, $\bbS^2_{0,2}$, and $(\PP^1)^2_{0,2}$} \label{fn2}
\end{figure}

The non-standard real toric Del Pezzo surfaces
$\bbS^2$, $\bbS^2_{1,0}$, $\bbS^2_{2,0}$, $\bbS^2_{0,2}$, and
$(\PP^1)^2_{0,2}$
%\footnote{The last surface is mentioned for
%completeness.}
can be associated with the following polygons
$\Del$, respectively (see Figure \ref{fn2}):
\begin{itemize}\item a square $\conv\{(0,0),(d,0),(0,d),(d,d)\}$,
$d\ge 1$; \item a pentagon
$\conv\{(0,0),(0,d),(d-d_1,d),(d,d-d_1),(d,0)\}$, $1\le d_1<d$;
\item a hexagon
$\conv\{(d_2,0),(0,d_2),(0,d),(d-d_1,d),(d,d-d_1),(d,0)\}$,
\mbox{$1\le d_1\le d_2<d$}; \item a hexagon
$\conv\{(0,0),(0,d-d_1),(d_1,d),(d,d),(d,d_1),(d_1,0)\}$, $1\le
d_1<d$; \item a hexagon
$\conv\{(0,0),(d_1-d_3,0),(d_1,d_3),(d_1,d_2),(d_3,d_2),(0,d_2-d_3)\}$,
\mbox{$1\le d_3<d_2\le d_1$}.
\end{itemize} For the first four surfaces, the conjugation acts in
the torus $(\C^*)^2$ by $\conj(x,y)=(\overline y,\overline x)$,
and acts in the tautological line bundle ${\cal L}_\Del$,
generated by monomials $x^iy^j$, $(i,j)\in\Del\cap\Z^2$, by
$\conj_*(a_{ij}x^iy^j)=\overline{a}_{ij}x^jy^i$, $(i,j)\in\Del$,
resembling the reflection of $\Del$ with respect to the bisectrix
${\cal B}$ of the positive quadrant\footnote{We shall denote by
${\cal B}$ the bisectrix of the positive quadrant both in the
plane of exponents of the polynomials in consideration, and in the
image-plane of the non-Archimedean valuation, no confusion will
arise. Furthermore, we shall denote by ${\cal B}_+$ and ${\cal
B}_-$ the halfplanes supported at ${\cal B}$ and lying
respectively above or below ${\cal B}$.}. For the fifth surface,
the action in $(\C^*)^2$ is $\conj(x,y)=(1/\overline x,1/\overline
y)$, and the action in ${\cal L}_\Del$ is
$\conj_*(a_{ij}x^iy^j)=\overline{a_{i,j}}x^{d_1-i}y^{d_2-j}$,
$(i,j)\in\Del$, resembling the reflection of $\Del$ with respect
to its center.

Observe that $-c_1({\cal L}_\Del)K_\Sig=|\partial\Del|$, the
lattice length of $\partial\Del$.

\subsection{Welschinger invariants associated with the totally real configurations of points}\label{secsp1}

\subsubsection{Admissible lattice paths and graphs}\label{secsp30}
For any lattice polygon $\del\subset\R^2$, symmetric with respect
to ${\cal B}$, denote by $(\partial\del)_\perp$ the union of the
sides of $\del$, orthogonal to ${\cal B}$, and denote by
$(\partial\del)_+$ the union of the sides of $\del\cal{\cal B}_+$
which are not orthogonal to ${\cal B}$. By
$\Tor((\partial\del)_\perp)$ (resp., $\Tor((\partial\del)_+)$) we
denote the union of the toric divisors $\Tor(\sig)$, where
$\sig\subset(\partial\del)_\perp$ (resp.,
$\sig\subset(\partial\del)_+$), in the toric surface $\Tor(\del)$,
associated with the polygon $\del$.

\begin{figure}
\setlength{\unitlength}{1cm}
\begin{picture}(13,8)(0,0)
\thinlines
\put(0,1){\vector(1,0){6}}\put(0,1){\vector(0,1){6}}\put(3,1){\line(0,1){1}}\put(4,1){\line(0,1){1}}
\put(1,5){\line(0,1){1}}\put(2,5){\line(0,1){1}}\put(0,5){\line(1,0){1}}\put(0,4){\line(1,0){1}}
\put(4,3){\line(1,0){1}}\put(4,2){\line(1,0){1}}
\thicklines\put(0,1){\line(0,1){2}}\put(0,3){\line(1,0){1}}\put(1,3){\line(0,1){2}}
\put(1,5){\line(1,0){2}}\put(3,5){\line(0,1){1}}\put(3,6){\line(1,0){2}}
\put(5,6){\line(0,-1){2}}\put(5,4){\line(-1,0){1}}\put(4,4){\line(0,-1){2}}
\put(4,2){\line(-1,0){2}}\put(0,1){\line(1,0){2}}\put(2,1){\line(0,1){1}}
\thinlines\dashline{0.2}(0,1)(5.5,6.5)
\dashline{0.2}(0,6)(3,6)\dashline{0.2}(5,1)(5,4)
\put(5.7,6.3){${\cal B}$}\put(3,0){\text{\rm
(a)}}\thinlines\put(1,1){\line(-1,1){1}}\put(2,1){\line(-1,1){2}}
\put(2,2){\line(-1,1){1}}\put(3,2){\line(-1,1){2}}\put(4,2){\line(-1,1){3}}
\put(4,3){\line(-1,1){2}}\put(4,4){\line(-1,1){1}}\put(5,4){\line(-1,1){2}}
\put(5,5){\line(-1,1){1}}\thinlines
\put(7,1){\vector(1,0){6}}\put(7,1){\vector(0,1){6}}\put(8,1){\line(-1,1){1}}
\put(9,1){\line(-1,1){2}}\put(9,2){\line(-1,1){1}}\put(10,3){\line(-1,1){1}}
\put(11,3){\line(-1,1){2}}\put(11,4){\line(-1,1){1}}\put(12,4){\line(-1,1){2}}
\put(12,5){\line(-1,1){1}}\thicklines
\put(7,1){\line(0,1){2}}\put(7,3){\line(1,0){4}}
\put(9,5){\line(1,0){1}}\put(10,5){\line(0,1){1}}\put(10,6){\line(1,0){2}}
\put(12,6){\line(0,-1){2}}\put(12,4){\line(-1,0){1}}\put(11,4){\line(0,-1){1}}
\put(7,1){\line(1,0){2}}\put(9,1){\line(0,1){4}}
\thinlines\dashline{0.2}(7,1)(12.5,6.5)
\dashline{0.2}(7,6)(10,6)\dashline{0.2}(12,1)(12,4)
\put(12.7,6.3){${\cal B}$}\put(10,0){\text{\rm (b)}}
\end{picture}
\caption{Lattice paths and subdivisions of $\Del$}\label{figsp2}
\end{figure}
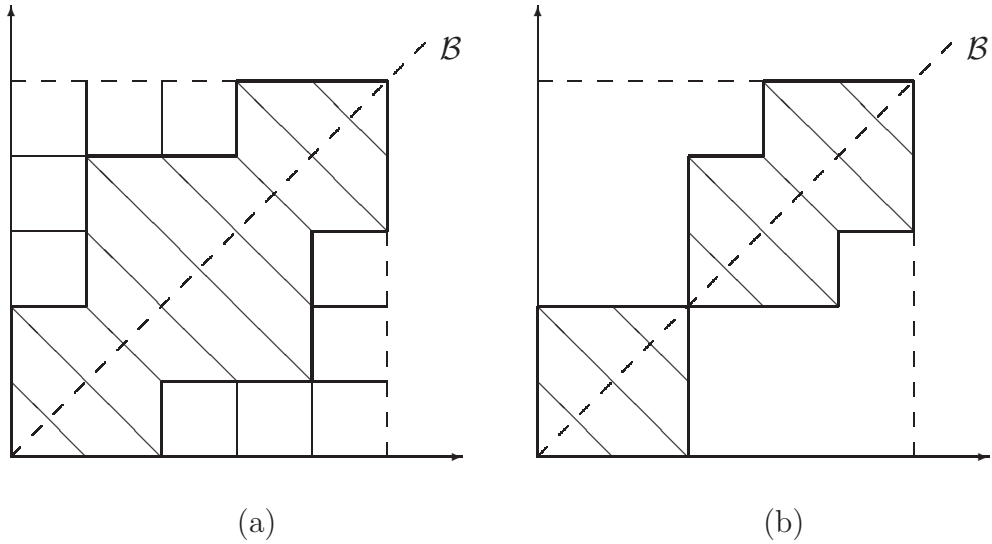

 Let $\Del$ be
one of the four polygons shown in Figure \ref{fn2}(a-d). The
integral points divide $(\partial\Del)_+$ into segments $s_i$,
$1\le i\le m:=|(\partial\Del)_+|$. An {\bf admissible lattice
path} in $\Del$ is a map $\gam:[0,m]\to\Del$ such that (see
example in Figure \ref{figsp2}(a))
\begin{itemize}\item image of $\gam$ lies in ${\cal
B}_+$;\item $\gam(0)$ and $\gam(m)$ are the two endpoints of
$(\partial\Del)_+$;\item the composition of the functional $x+y$
with $\gam$ is a strongly increasing function;\item
$\gam(i)\in\Z^2$ , and $\gam\big|_{[i,i+1]}$ is linear as
$i\in\Z$;\item there is a permutation $\tau\in S_{m-1}$ such that
$\gam([i-1,i])$ is a translate of the segment $s_{\tau(i)}$, $1\le
i\le m$, of $(\partial\Del)_+$;\item $\gam([0,m])\cap{\cal
B}=(\partial\Del)_+\cap{\cal B}$.\end{itemize} The image of $\gam$
completely determines the map, and we shall denote them by the
same symbol $\gam$.

Denote by $\sig_i$, $0\le i\le m$, the segments, joining the
integral points of $\gam$ with their symmetric images (see Figure
\ref{figsp2}(a)).

A {\bf $\gam$-admissible graph} $G$ is defined as follows. First,
we describe some subgraph $G'$. The connected components of $G'$
are lattice segments (or points)
\mbox{$G'_j=[(a_j,j),(b_j,j)]\subset\R^2$}, $1\le j\le
n:=|\partial\Del|-m-1$, with positive integer weights $w(G'_j)$
such that \begin{itemize}\item $0\le a_j\le b_j\le m$ for all
$j=1,...,n$; \item $a_j\le a_{j+1}$, and in addition $b_j\le
b_{j+1}$ if $a_j=a_{j+1}$ as $j=1,...,n-1$; \item for all
$i=0,...,m$, \begin{equation}\sum_{(i,j)\in G'_j}w(G'_j)=|\sig_i|\
;\label{esp21}\end{equation} \item if $a_j=0$ or $b_j=m$ then
$w(G'_j)=1$. \end{itemize}

We then introduce new vertices $\varphi_i$, $i=1,...,m$, of the
graph $G$ and the new arcs, joining any vertex $\varphi_i$ with
the endpoint $(i-1,j)$ of any component $G'_j$ satisfying
$b_j=i+1$, and with the endpoint $(i,j)$ of any component $G'_j$
satisfying $a_j=i$. Our final requirement is that the obtained
graph $G$ is a tree.

A {\bf marking} of a $\gam$-admissible graph $G$ is a vector
$\overline s=(s_1,...,s_n)\in\Z^n$ such that $a_j\le s_j\le b_j$,
$j=1,...,n$, subject to the following restriction:
$$s_j\le s_{j+1}\quad\text{\rm as far as}\quad a_j=a_{j+1},\
b_j=b_{j+1}\ .$$

At last we define the {\bf Welschinger number} $W(\gam,G,\overline
s)=0$ if at leat one weight $w(G'_i)$ is even, and otherwise
\begin{equation}W(\gam,G,\overline
s)=2^v\prod_{k=0}^mn_k!\left(\prod_{\renewcommand{\arraystretch}{0.6}
\begin{array}{c}
\scriptstyle{0\le a\le b\le m}\\
\scriptstyle{c=1,3,5,...}
\end{array}}
n_{k,a,b,c}!\right)^{-1}\ ,\label{espn60}\end{equation} where $v$
is the total valency of those vertices $\varphi_i$ of $G$, for
which $|\sig_i|=|\sig_{i-1}|$, and
$$n_k:=\#\{j\ :\ s_j=k,\ 1\le j\le n\},\quad n_{k,a,b,c}=\#\{j\ :\
s_j=k,\ a_j=a,\ b_j=b,\ w(G'_j)=c\},$$ $$k=0,...,m,\quad 0\le a\le
k\le b\le m,\quad c=1,3,5,...$$

\subsubsection{The formula}\label{secsp3}

\begin{theorem}\label{tsp1}
In the notation of sections \ref{sec1} and \ref{secsp30}, if
$\Sig=\bbS^2$, $\bbS^2_{1,0}$, $\bbS^2_{2,0}$, or $\bbS^2_{0,2}$,
then
\begin{equation}W_0(\Sig,{\cal L}_\Del)=\sum W(\gam,G,\overline s)\
,\label{esp22}\end{equation} where the sum ranges over all
admissible lattice paths $\gam$, all $\gam$-admissible graphs $G$,
and all markings $\overline s$ of $G$.
\end{theorem}

It is an easy exercise to show that there always exist an
admissible lattice path and a corresponding admissible graph, and
hence

\begin{corollary}\label{csp2} In the above notation,
for any surface $\Sig=\bbS^2$, $\bbS^2_{1,0}$, $\bbS^2_{2,0}$, or
$\bbS^2_{0,2}$, and any line bundle ${\cal L}_\Del$, the
Welschinger invariant $W_0(\Sig,{\cal L}_\Del)$ is positive, and
through any $-c_1({\cal L}_\Del)K_\Sig-1$ generic real points on
$\Sig$ there passes at least one real rational curve $D\in|{\cal
L}_\Del|$.
\end{corollary}

\subsubsection{Examples}\label{secspn2}

\begin{figure}
\setlength{\unitlength}{1cm}
\begin{picture}(14,4)(0,0)
\thinlines
\put(0.5,1){\vector(1,0){3}}\put(0.5,1){\vector(0,1){3}}
\put(4,1){\vector(1,0){3}}\put(4,1){\vector(0,1){3}}
\put(7.5,1){\vector(1,0){3}}\put(7.5,1){\vector(0,1){3}}
\put(11,1){\vector(1,0){3}}\put(11,1){\vector(0,1){3}}
\put(2.5,1){\line(0,1){2}}\put(6,1){\line(0,1){1}}\put(6,2){\line(-1,1){1}}\put(7.5,2){\line(1,-1){1}}
\put(9.5,1){\line(0,1){1}}\put(8.5,3){\line(1,-1){1}}
\thicklines\put(12,1){\line(1,1){1}}\put(13,2){\line(0,1){1}}
\put(0.5,1){\line(0,1){2}}\put(0.5,3){\line(1,0){2}}
\put(4,1){\line(0,1){2}}\put(4,3){\line(1,0){1}}
\put(7.5,2){\line(0,1){1}}\put(7.5,3){\line(1,0){1}}
\put(11,1){\line(0,1){1}}\put(12,3){\line(1,0){1}}\put(12,3){\line(-1,-1){1}}
\thinlines\dashline{0.2}(0.5,2)(1.5,1)
\dashline{0.2}(0.5,3)(2.5,1)\dashline{0.2}(1.5,3)(2.5,2)
\dashline{0.2}(4,2)(5,1)\dashline{0.2}(4,3)(6,1)\dashline{0.2}(7.5,3)(9.5,1)
\dashline{0.2}(11,2)(12,1)\dashline{0.2}(12,3)(13,2)
\dashline{0.1}(11,3)(12,3)\dashline{0.1}(13,2)(13,1)\put(0.2,2.9){$2$}\put(2.4,0.6){$2$}
\put(3.7,2.9){$2$}\put(5.9,0.6){$2$}\put(7.2,2.9){$2$}\put(9.4,0.6){$2$}
\put(10.7,2.9){$2$}\put(12.9,0.6){$2$}
\end{picture}
\caption{Linear systems with an elliptic general
curve}\label{figsp4}
\end{figure}
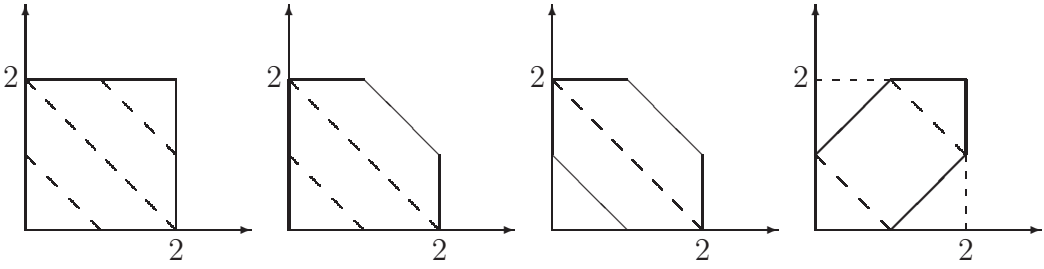

{\it (A) Linear systems with an elliptic general member.} Let
$\Sig=\bbS^2$, $\bbS^2_{1,0}$, $\bbS^2_{2,0}$, or $\bbS^2_{0,2}$,
and let $\Del$ be a respective associated lattice polygon as shown
in Figure \ref{fn2} and such that a general curve in $|{\cal
L}_\Del|$ is elliptic. Then $\Del$ is as depicted in Figure
\ref{figsp4}.

The Welschinger invariant $W_0(\Sig,{\cal L}_\Del)$ can be
computed by counting rational curves in the pencil of real curves
in $|{\cal L}_\Del|$ passing through $(-c_1({\cal
L}_\Del)K_\Sig-1)$ generic real points. Integrating along the
pencil with respect to the Euler characteristic and noticing that
the curves in the pencil have one more real base point, and
$\chi(\R D)=1$ or $-1$ according as $D\in|{\cal L}_\Del|$ is a
real rational curve with a solitary or a non-solitary node, we
obtain (cf. with the case of plane cubics \cite{IKS2}, section
3.1) $$W_0(\Sig,{\cal L}_\Del)=-c_1({\cal
L}_\Del)K_\Sig-\chi(\Sig(\R))\ ,$$ which equals $6$, $6$, $6$, or
$4$ as $\Sig=\bbS^2$, $\bbS^2_{1,0}$, $\bbS^2_{2,0}$, or
$\bbS^2_{0,2}$, respectively.

In turn in Theorem \ref{tsp1} we have a unique admissible path
$\gam$ (fat line in Figure \ref{figsp4}) and a unique subdivision
of $\Del$ (dashes in Figure \ref{figsp4}). The subgraphs $G'$ of
the $\gam$-admissible graphs, their markings and Welschinger
numbers are shown in Figure \ref{figsp5} (the weight of any
component of $G'$ is here $1$). The result, of course, coincides
with the aforementioned one.

\begin{figure}
\setlength{\unitlength}{1cm}
\begin{picture}(14,10)(0,0)
\thinlines
\put(5.5,1.5){\vector(1,0){2}}\put(0,6.5){\vector(1,0){3}}
\put(3.5,6.5){\vector(1,0){3}}\put(7,6.5){\vector(1,0){3}}
\put(10.5,6.5){\vector(1,0){3}}\dashline{0.1}(5.5,2.5)(7.5,2.5)
\dashline{0.1}(5.5,3.5)(7.5,3.5)\dashline{0.1}(6,1.5)(6,3.7)
\dashline{0.1}(7,1.5)(7,3.7)\put(5.91,2.4){$\bullet$}\put(6.91,3.4){$\bullet$}
\dashline{0.1}(0,7.5)(3,7.5)
\dashline{0.1}(0,8.5)(3,8.5)\dashline{0.1}(0,9.5)(3,9.5)
\dashline{0.1}(0.5,6.5)(0.5,9.7)\dashline{0.1}(1.5,6.5)(1.5,9.7)\dashline{0.1}(2.5,6.5)(2.5,9.7)
\dashline{0.1}(3.5,7.5)(6.5,7.5)
\dashline{0.1}(3.5,8.5)(6.5,8.5)\dashline{0.1}(3.5,9.5)(6.5,9.5)
\dashline{0.1}(4,6.5)(4,9.7)\dashline{0.1}(5,6.5)(5,9.7)\dashline{0.1}(6,6.5)(6,9.7)
\dashline{0.1}(7,7.5)(10,7.5)
\dashline{0.1}(7,8.5)(10,8.5)\dashline{0.1}(7,9.5)(10,9.5)
\dashline{0.1}(7.5,6.5)(7.5,9.7)\dashline{0.1}(8.5,6.5)(8.5,9.7)\dashline{0.1}(9.5,6.5)(9.5,9.7)
\dashline{0.1}(10.5,7.5)(13.5,7.5)
\dashline{0.1}(10.5,8.5)(13.5,8.5)\dashline{0.1}(10.5,9.5)(13.5,9.5)
\dashline{0.1}(11,6.5)(11,9.7)\dashline{0.1}(12,6.5)(12,9.7)\dashline{0.1}(13,6.5)(13,9.7)
\put(0.41,7.4){$\bullet$}\put(1.41,8.4){$\bullet$}\put(1.41,9.4){$\bullet$}
\put(3.91,7.4){$\bullet$}\put(4.91,8.4){$\bullet$}\put(5.91,9.4){$\bullet$}
\put(7.41,7.4){$\bullet$}\put(8.41,8.4){$\bullet$}\put(9.41,9.4){$\bullet$}
\put(11.91,7.4){$\bullet$}\put(11.91,8.4){$\bullet$}\put(12.91,9.4){$\bullet$}
\thicklines\put(1.5,9.5){\line(1,0){1}}\put(5,9.5){\line(1,0){1}}\put(7.5,7.5){\line(1,0){1}}
\put(11,7.5){\line(1,0){1}}
\put(1,5.5){$W=2$}\put(4.5,5.5){$W=1$}\put(8,5.5){$W=1$}\put(11.5,5.5){$W=2$}
\put(5,4.5){$\text{\rm (a)}\quad\Sig=\bbS^2,\ \bbS^2_{1.0},\
\text{\rm or}\
\bbS^2_{2,0}$}\put(8.5,2.5){$W=4$}\put(6,0.5){$\text{\rm
(b)}\quad\Sig=\bbS^2_{0,2}$}
\end{picture}
\caption{Admissible graphs and markings, I}\label{figsp5}
\end{figure}

{\it (B) Linear systems of digonal curves.} We illustrate Theorem
\ref{tsp1} by two more examples, where one can easily write down a
closed formula for the Welschinger invariant (a similar
computation has been performed for digonal curves on $(\PP^1)^2$
\cite{IKS2}, section 3.1). Namely, we consider the surfaces
$\Sig=\bbS^2_{2,0}$ and $\bbS^2_{0,2}$ and the linear systems
associated with the polygons $\Del$ shown in Figure
\ref{figsp6}(a,b), respectively.

\begin{figure}
\setlength{\unitlength}{0.95cm}
\begin{picture}(14,22)(0,0)
\thinlines
\put(1,17.5){\vector(1,0){5}}\put(8,17.5){\vector(1,0){7}}
\put(1,10){\vector(1,0){5}}\put(8,10){\vector(1,0){5}}
\put(1,17.5){\vector(0,1){5}}\put(8,17.5){\vector(0,1){5}}
\put(1,10){\vector(0,1){5}}\put(8,10){\vector(0,1){5}}
\put(1,1){\vector(1,0){5}}\put(9,1){\vector(1,0){3}}
\dashline{0.1}(1,21.5)(4,21.5) \dashline{0.1}(5,17.5)(5,20.5)
\dashline{0.1}(1,2)(6,2) \dashline{0.1}(1,3)(6,3)
\put(1.41,1.9){$\bullet$}\put(2.41,2.9){$\bullet$}\put(3.41,3.9){$\bullet$}
\put(4.41,4.9){$\bullet$} \dashline{0.1}(1,4)(6,4)
\dashline{0.1}(1,5)(6,5)\dashline{0.1}(1.5,1)(1.5,5.2)
\dashline{0.1}(2.5,1)(2.5,5.2)\dashline{0.1}(3.5,1)(3.5,5.2)\dashline{0.1}(4.5,1)(4.5,5.2)
\dashline{0.1}(9,2)(12,2)
\dashline{0.1}(9,3)(12,3)\dashline{0.1}(9,4)(12,4)
\dashline{0.1}(9,5)(12,5)\dashline{0.1}(9,6)(12,6)\dashline{0.1}(9,7)(12,7)
\dashline{0.1}(9,8)(12,8)\dashline{0.1}(8,18.5)(12,18.5)
\dashline{0.1}(9.5,1)(9.5,8.2)\dashline{0.1}(10.5,1)(10.5,8.2)
\dashline{0.1}(11.5,1)(11.5,8.2)\put(9.41,1.9){$\bullet$}\put(9.41,2.9){$\bullet$}
\put(10.41,3.9){$\bullet$}\put(10.41,4.9){$\bullet$}\put(10.41,5.9){$\bullet$}
\put(11.41,6.9){$\bullet$}\put(11.41,7.9){$\bullet$}
\put(4,21.5){\line(1,0){1}}\put(1,18.5){\line(1,1){3}}\put(2,17.5){\line(1,1){3}}
\put(5,20.5){\line(0,1){1}}\put(8,20.5){\line(1,-1){3}}\put(9,21.5){\line(1,-1){3}}
\put(8,21.5){\line(1,0){1}}\put(12,17.5){\line(0,1){1}}\put(1,13){\line(1,-1){3}}
\put(2,14){\line(1,-1){3}}\put(5,10){\line(0,1){1}}\put(8,13){\line(1,-1){3}}
\put(9,14){\line(1,-1){3}}\put(9,14){\line(-1,0){1}}\put(12,10){\line(0,1){1}}
\put(0.6,18.4){$1$}\put(0.6,21.4){$d$}\put(4.9,17){$d$}\put(1.9,17){$1$}\put(7,20.4){$d-1$}
\put(7.6,21.4){$d$}\put(11.9,17){$d$}\put(7.6,18.4){$1$}\put(3.4,16){$\text{\rm
(a)}$}\put(10.4,16){$\text{\rm (b)}$}\put(3.4,9){$\text{\rm
(c)}$}\put(10.4,9){$\text{\rm (d)}$}\put(3.4,0){$\text{\rm
(e)}$}\put(10.4,0){$\text{\rm
(f)}$}\put(1.2,13.5){$\gam_1$}\put(9.1,12.6){$\gam_2$}
\thicklines\put(1,13){\line(0,1){1}}\put(1,14){\line(1,0){1}}\put(8,13){\line(1,0){1}}
\put(9,13){\line(0,1){1}}\put(9.5,3){\line(1,0){1}}\put(9.5,4){\line(1,0){1}}
\put(10.5,6){\line(1,0){1}}\put(10.5,7){\line(1,0){1}}
\end{picture}
\caption{Admissible paths, graphs, and markings, II}\label{figsp6}
\end{figure}

In the case $\Sig=\bbS^2_{0,2}$ (see Figure \ref{figsp6}(b)) we
have a unique admissible path $\gam$ going just along
$(\partial\Del)_+$, a unique $\gam$-admissible graph $G$, and a
unique marking (see Figure \ref{figsp6}(e)). Hence we obtain
$W_0(\bbS^2_{0,2},{\cal L}_\Del)=4^{d-1}$, $d$ being the length of
projection of $\Del$ on a coordinate axis.

In the case $\Sig=\bbS^2_{2,0}$, $d>2$, there are two admissible
lattice paths $\gam_1,\gam_2$ (shown by fat lines in Figure
\ref{figsp6}(c,d)). The subgraph $G'$ of an admissible graph $G$,
and a marking $\overline s$ should look as shown in Figure
\ref{figsp6}(f), where we denote by $k$ (resp., $l$) the number of
components $[(1,j),(2,j)]$ (resp. $[(2,j),(3,j)]$), and $k_1$
(resp. $l_1$) is the number of the marked points $(2,j)$ on
components $[(1,j),(2,j)]$ (resp. $[(2,j),(3,j)]$), and where
$(k,l,k_1,l_1)$ run over the sets $J(\gam_1)$ and $J(\gam_2)$
defined by
$$0\le k_1\le k\le d-1,\quad 0\le l_1\le l\le d-1,\quad
d-k-l=2a+1,\ a\ge\begin{cases}1,\ &\text{\rm if}\ \gam=\gam_1,\\
2,\ &\text{\rm if}\ \gam=\gam_2\end{cases}$$ Here the weights of
all the components of $G'$ are equal to $1$, except for the
one-point component on the middle vertical line, whose weight is
$2a+1$ or $2a-1$ according as $\gam=\gam_1$ or $\gam_2$. Thus, we
obtain
$$W_0(\bbS^2_{2,0},{\cal L}_\Del)=\sum_{i=1}^2\sum_{(k,l,k_1,l_1)
\in
J(\gam_i)}\frac{(d-1-k_1)!(d-1-l_1)!(k_1+l_1+1)!}{(d-1-k)!(d-1-l)!(k-k_1)!(l-l_1)!k_1!l_1!}\
.$$

\subsection{Welschinger invariants associated with configurations containing imaginary points}

Let $\Sig$ be one of the surfaces $\bbS^2$, $\bbS^2_{1,0}$,
$\bbS^2_{2,0}$, or $\bbS^2_{0,2}$, and let $\Del$ be a respective
lattice polygon shown in Figure \ref{fn2}(a-d). Fix positive
integers $r',r''$ satisfying (\ref{enn36}) in the form
\begin{equation}r'+2r''=|\partial\Del|-1\
.\label{espn1}\end{equation}

\subsubsection{The case of $\Sig=\bbS^2$, $\bbS^2_{1,0}$,
or $\bbS^2_{2,0}$}\label{secspn1} First, we introduce a splitting
${\cal R}$ of $r'$ and $r''$ into nonnegative integer summands
\begin{equation}r'=r'_1+r'_2,\quad r''=r''_1+r''_{2,1}+r''_{2,2}\
,\label{espn3}\end{equation} so that (cf. (\ref{espn1}),
(\ref{espn3}))
\begin{equation}r'_1+2r''_1+r''_{2,1}=|\partial\Del|-|\partial\Del|_+-1,\quad
r'_2+r''_{2,1}+2r''_{2,2}=|\Del|_+\ .\label{espn2}\end{equation}
Notice here that, for $r''=0$, the splitting ${\cal R}$ turns into
\begin{equation}r'=|\partial\Del|-1,\quad
r'_1=|\partial\Del|-|\partial\Del|_+-1,\quad r'_2=|\Del|_+\
.\label{espn53}\end{equation}

Take a lattice path $\gam$ in $\Del$, admissible in the sense of
section \ref{secsp30}. We then pick integral points $v_i$,
$i=0,1,...,\widetilde m:=r'_2+r''_{2,1}+r''_{2,2}$, on $\gam$ such
that
\begin{enumerate}\item[($\gam${\bf1})] they are ordered by the growing sum of
coordinates; \item[($\gam${\bf2})] $v_0$ and $v_{\widetilde m}$
are the endpoints;
\item[($\gam${\bf3})] the lattice length of the part $\gam_i$ of $\gam$
between the points $v_{i-1},v_i$, $1\le i\le\widetilde m$, is $1$
or $2$;
\item[($\gam${\bf4})] $|\gam_i|=1$ for all $i>r''_{2,1}+r''_{2,2}$;\item[($\gam${\bf5})]
each break point of $\gam$, where the path turns in positive
direction (i.e, left), is one of $v_i$, $1\le i\le\widetilde m$.
\end{enumerate} Denote by $\sig_i$ the segment, joining $v_i$ with
its symmetric (with respect to ${\cal B}$) image,
$i=0,...,\widetilde m$. Observe that $\gam$ defines a lattice
subdivision $S$ of $\Del$, consisting of the polygons in the part
of $\Del$ between $\gam$ and its symmetric image $\hat\gam$, cut
out by the segments $\sig_i$, $i=1,...,\widetilde m-1$, and the
rectangles with vertical and horizontal sides in the remaining
part of $\Del$ (see Figure \ref{figsp2}(a)).

Then we construct a $(\gam,{\cal R})$-admissible graph $G$,
starting with an auxiliary subgraph $G'$, whose connected
components are segments or points $G'_j=[(a_j,j),(b_j,j)]$ lying
on the lines $y=j$ with
$j=1,...,n:=|\partial\Del|-|\partial\Del|_+-1=r'_1+2r''_1+r''_{2,1}$,
respectively, equipped with positive integer weights $w(G'_j)$,
and satisfying the following conditions:
\begin{enumerate}\item[({\bf G1})] for all $j=1,...,n$,
$$0\le a_j\le b_j\le\widetilde m\ ;$$\item[({\bf G2})] for all
$j=1,...,r''_1$,
$$a_{2j-1}=a_{2j},\quad b_{2j-1}=b_{2j},\quad w(G'_{2j-1})=w(G'_{2j})\
;$$\item[({\bf G3})] $w(G'_j)=1$ as far as $a_j=0$ or
$b_j=\widetilde m$;
\item[({\bf G4})] $a_j\le a_{j+1}$ and, if $a_j=a_{j+1}$ then $b_j\le
b_{j+1}$, as far as either $1\le j<2r''_1$, or
$2r''_{2,1}<j<2r''_1+r''_{2,1}$, or $2r''_1+r''_{2,1}<j<
n$;\item[({\bf G5})] if $|\gam_i|=2$ for some $1\le i\le\widetilde
m$, then
$$a_j\ne i+1\ \text{and}\ b_j\ne i\quad\text{for all}\ 2r''_1<j\le
n\ ;$$\item[({\bf G6})] $a_j\le r''_{2,1}+r''_{2,2}$ for
$2r''_1<j\le 2r''_1+r''_{2,1}$; \item[({\bf G7})] for any
$i=0,...,\widetilde m$, relation (\ref{esp21}) holds true.
\end{enumerate} Next we introduce some more vertices of $G$, taking
one vertex $\varphi_i$ for all $i=1,...,\widetilde m$ such that
$|\gam_i|=1$, and taking two vertices
$\varphi_{i,1},\varphi_{i,2}$ for all $i=1,...,\widetilde m$ such
that $|\gam_i|=2$, and, finally, we define additional arcs of $G$
as follows:
\begin{enumerate}\item[({\bf G8})] any vertex $\varphi_i$, $1\le
i\le\widetilde m$, is joined by arcs with each endpoint $(a_j,j)$,
$a_j=i$, of a component of $G'$, and with each endpoint $(b_j,j)$,
$b_j=i-1$, of a component of $G'$; \item[({\bf G9})] any vertex
$\varphi_{i,1}$ (resp., $\varphi_{i,2}$) is joined by arcs with
each endpoint $(a_{2j-1},2j-1)$ (resp., $(a_{2j},2j)$),
$a_{2j-1}=i$, $1\le j\le r''_1$, of a component of $G'$, and with
each endpoint $(b_{2j-1},2j-1)$ (resp., $(b_{2j},2j)$),
$b_{2j}=i-1$, $1\le j\le r''_1$, of a component of $G'$;
\item[({\bf G10})]the minimal subgraph $G''$ of $G$, containing the components
$G'_j$, $2r''_1<j\le2r''_1+r''_{2,1}$, and all the vertices
$\varphi_i$, $1\le i\le r''_{2,1}+r''_{2,2}$, is a forest, each
component of $G''$ contains at least one univalent vertex of type
$\varphi_i$, furthermore, any component of $G''$ can be oriented
so that from each vertex of type $\varphi_i$ emanates precisely
one oriented arc, and no arc emanates from a vertex of type
$(0,j)$ or $(\widetilde m,j)$;
\item[({\bf G11})] $G$ is a tree.
\end{enumerate}

A marking of a $(\gam,{\cal R})$-admissible graph $G$ is a pair of
integer vectors $\overline s'=(s'_1,...,s'_{r'_1})$ and $\overline
s''=(s''_1,...,s''_{r''_1})$ such that
\begin{enumerate}\item[({\bf s1})] $a_{j+2r''_1+r''_{2,1}}\le s'_j\le b_{j+2r''_1+r''_{2,1}}$ as
$j=1,...,r'_1$, and $a_{2j}\le s''_j\le b_{2j}$ as
$j=1,...,r''_1$;\item[({\bf s2})] $s'_j\le s'_{j+1}$ as far as
$a_{j+2r''_1+r''_{2,1}}=a_{j+2r''_1+r''_{2,1}+1},\
b_{j+2r''_1+r''_{2,1}}=b_{j+2r''_1+r''_{2,1}+1}$, where $1\le
j<r'_1$;
\item[({\bf s3})] $s''_j\le s''_{j+1}$ as far as $a_{2j}=a_{2j+1},\
b_{2j}=b_{2j+1}$, where $1\le j<r''_1$;
\item[({\bf s4})] $s'_j\ge s''_k$ for all $j=1,...,r'_1$ and $k=1,...,r''_1$.
\end{enumerate} Observe that condition ({\bf s4}) imposes an extra
restriction to the components of $G'$ and splitting ${\cal R}$.

Then we define the Welschinger number as $W({\cal
R},\gam,G,\overline s',\overline s'')=0$ if at least one weight
$w(G'_i)$, $2r''_1<i\le n$, is even, and otherwise as
\begin{equation}W({\cal R},\gam,G,\overline s',\overline
s'')=(-1)^a\cdot
2^{b+c}\cdot\prod_{j=1}^{r''_1}(w(G'_{2j}))^2\cdot\prod_{j=2r''_1+1}^{2r''_1+r''_{2,1}}
w(G'_j)\cdot\prod_{i=0}^{\widetilde m}
\left(n'_i!n''_i!\alp_i^{-1}\bet_i^{-1}\right)\
,\label{espn61}\end{equation}
$$\alp_i=\prod_{\renewcommand{\arraystretch}{0.6}
\begin{array}{c}
\scriptstyle{0\le d\le e\le\widetilde m}\\
\scriptstyle{f=1,3,5,...}
\end{array}}
n'_{i,d,e,f}!,\quad
\bet_i=\prod_{\renewcommand{\arraystretch}{0.6}
\begin{array}{c}
\scriptstyle{0\le d\le e\le\widetilde m}\\
\scriptstyle{f=1,2,3,...}
\end{array}}
n''_{i,d,e,f}!\ ,$$ where \begin{itemize}\item
$a=(a_1-a'_1)+...+(a_{\widetilde m}-a'_{\widetilde m})$ with the
following summands: if $|\gam_k|=1$, then $a_k=a'_k=0$, if
$|\gam_k|=2$, then $a_k$ is the length of the segment $[p,\hat
p]$, $p$ being the intermediate integral point of $\gam$, and
$a'_k$ is the total weight of the components of $G'$, crossing the
vertical line $x=k-1/2$;
\item $b$ is the total valency of all the vertices
$\varphi_{i,1}$ of $G$; \item $c$ is the number of hexagons in the
subdivision $S$;
\item $n'_i=\#\{j\ :\ s'_j=i,\ 1\le j\le r'_1\}$, $n''_i=\#\{j\ :\
s''_j=i,\ 1\le j\le r''_1\}$;\item $n'_{i,d,e,f}=\#\{j\ :\
s'_j=i,\ a_j=d,\ b_j=e,\ w(G'_j)=f\}$;
\item $n''_{i,d,e,f}=\#\{j\ :\ s''_j=i,\ a_{2j}=d,\
b_{2j}=e,\ w(G'_{2j})=f\}$.
\end{itemize}

\begin{theorem}\label{tsp3}
In the notation of section \ref{secspn1}, if $\Sig=\bbS^2$,
$\bbS^2_{1,0}$, or $\bbS^2_{2,0}$, and the positive integers
$r',r''$ satisfy (\ref{espn1}), then
\begin{equation}W_{r''}(\Sig,{\cal L}_\Del)=\sum W({\cal R},\gam,G,\overline s',\overline s'')\
,\label{espn4}\end{equation} where the sum ranges over splittings
${\cal R}$ of $r',r''$, satisfying (\ref{espn2}), all admissible
lattice paths $\gam$, all $\gam$-admissible graphs $G$, and all
markings $\overline s',\overline s''$ of $G$, subject to
conditions, specified in section \ref{secspn1}.
\end{theorem}

As example we consider the linear system of curves of bi-degree
$(2,2)$ on $\bbS^2$. Here $r'+2r''=7$, and the integration with
respect to Euler characteristic (cf. section \ref{secspn2}, case
(A)) gives $W_{r''}(\bbS^2,(2,2))=6-2r''$, $r''=1,2,3$. In Figure
\ref{figsp10} we demonstrate how to obtain these answers from
formula (\ref{espn4}), listing admissible paths $\gam$,
subdivisions $S$ of $\Del$, and graphs $G'$ (here the marked
points are denoted by bullets and the non-marked one-point
components of $G'$ are denoted by circles).

\begin{figure}
\setlength{\unitlength}{0.6cm}
\begin{picture}(20,35)(0,-1)
\thinlines \put(4,30){\vector(1,0){4}}\put(4,30){\vector(0,1){4}}
\put(9,30){\vector(1,0){4}}\put(9,30){\vector(0,1){4}}\put(14,30){\vector(1,0){4}}\put(14,30){\vector(0,1){4}}
\put(3,21.5){\vector(1,0){4}}\put(3,21.5){\vector(0,1){4}}\put(8,21.5){\vector(1,0){4}}\put(8,21.5){\vector(0,1){4}}
\put(13,21.5){\vector(1,0){4}}\put(13,21.5){\vector(0,1){4}}\put(18,21.5){\vector(1,0){4}}\put(18,21.5){\vector(0,1){4}}
\put(3,14.5){\vector(1,0){3}}\put(3,14.5){\vector(0,1){4}}
\put(4,7){\vector(1,0){4}}\put(4,7){\vector(0,1){4}}\put(10,7){\vector(1,0){4}}\put(10,7){\vector(0,1){4}}
\put(4,1){\vector(1,0){3}}\put(4,1){\vector(0,1){4}}
\dashline{0.1}(4,31)(7.5,31)\dashline{0.1}(4,32)(7.5,32)\dashline{0.1}(4,33)(7.5,33)
\dashline{0.1}(9,31)(12.5,31)\dashline{0.1}(9,32)(12.5,32)\dashline{0.1}(9,33)(12.5,33)
\dashline{0.1}(14,31)(17.5,31)\dashline{0.1}(14,32)(17.5,32)\dashline{0.1}(14,33)(17.5,33)
\dashline{0.1}(5,30)(5,33.5)\dashline{0.1}(6,30)(6,33.5)\dashline{0.1}(7,30)(7,33.5)
\dashline{0.1}(10,30)(10,33.5)\dashline{0.1}(11,30)(11,33.5)\dashline{0.1}(12,30)(12,33.5)
\dashline{0.1}(15,30)(15,33.5)\dashline{0.1}(16,30)(16,33.5)\dashline{0.1}(17,30)(17,33.5)
\dashline{0.1}(4,21.5)(4,25)\dashline{0.1}(5,21.5)(5,25)\dashline{0.1}(6,21.5)(6,25)
\dashline{0.1}(9,21.5)(9,25)\dashline{0.1}(10,21.5)(10,25)\dashline{0.1}(11,21.5)(11,25)
\dashline{0.1}(14,21.5)(14,25)\dashline{0.1}(15,21.5)(15,25)\dashline{0.1}(16,21.5)(16,25)
\dashline{0.1}(19,21.5)(19,25)\dashline{0.1}(20,21.5)(20,25)\dashline{0.1}(21,21.5)(21,25)
\dashline{0.1}(3,22.5)(6.5,22.5)\dashline{0.1}(3,23.5)(6.5,23.5)\dashline{0.1}(3,24.5)(6.5,24.5)
\dashline{0.1}(8,22.5)(11.5,22.5)\dashline{0.1}(8,23.5)(11.5,23.5)\dashline{0.1}(8,24.5)(11.5,24.5)
\dashline{0.1}(13,22.5)(16.5,22.5)\dashline{0.1}(13,23.5)(16.5,23.5)\dashline{0.1}(13,24.5)(16.5,24.5)
\dashline{0.1}(18,22.5)(21.5,22.5)\dashline{0.1}(18,23.5)(21.5,23.5)\dashline{0.1}(18,24.5)(21.5,24.5)
\dashline{0.1}(4,14.5)(4,18)\dashline{0.1}(5,14.5)(5,18)
\dashline{0.1}(3,15.5)(5.5,15.5)\dashline{0.1}(3,16.5)(5.5,16.5)\dashline{0.1}(3,17.5)(5.5,17.5)
\dashline{0.1}(7,7)(7,10.5)\dashline{0.1}(5,7)(5,10.5)\dashline{0.1}(6,7)(6,10.5)
\dashline{0.1}(11,7)(11,10.5)\dashline{0.1}(12,7)(12,10.5)\dashline{0.1}(13,7)(13,10.5)
\dashline{0.1}(4,8)(7.5,8)\dashline{0.1}(4,9)(7.5,9)\dashline{0.1}(4,10)(7.5,10)
\dashline{0.1}(10,8)(13.5,8)\dashline{0.1}(10,9)(13.5,9)\dashline{0.1}(10,10)(13.5,10)
\dashline{0.1}(5,1)(5,4.5)\dashline{0.1}(6,1)(6,4.5)
\dashline{0.1}(4,2)(6.5,2)\dashline{0.1}(4,3)(6.5,3)\dashline{0.1}(4,4)(6.5,4)
\dashline{0.1}(0,31)(1,30)\dashline{0.1}(0,32)(2,30)\dashline{0.1}(1,32)(2,31)
\dashline{0.1}(0,22.5)(1,21.5)\dashline{0.1}(0,23.5)(2,21.5)\dashline{0.1}(1,23.5)(2,22.5)
\dashline{0.1}(0,16.5)(2,14.5)\dashline{0.1}(1,16.5)(2,15.5)
\dashline{0.1}(0,8)(1,7)\dashline{0.1}(0,9)(2,7)\dashline{0.1}(1,9)(2,8)
\dashline{0.1}(0,3)(2,1)\dashline{0.1}(1,3)(2,2)
\put(0,30){\line(1,0){2}}\put(2,30){\line(0,1){2}}
\put(0,21.5){\line(1,0){2}}\put(2,21.5){\line(0,1){2}}
\put(0,14.5){\line(1,0){2}}\put(2,14.5){\line(0,1){2}}
\put(0,7){\line(1,0){2}}\put(2,7){\line(0,1){2}}
\put(0,1){\line(1,0){2}}\put(2,1){\line(0,1){2}} \thicklines
\put(0,30){\line(0,1){2}}\put(0,32){\line(1,0){2}}
\put(0,21.5){\line(0,1){2}}\put(0,23.5){\line(1,0){2}}
\put(0,14.5){\line(0,1){2}}\put(0,16.5){\line(1,0){2}}
\put(0,7){\line(0,1){2}}\put(0,9){\line(1,0){2}}
\put(0,1){\line(0,1){2}}\put(0,3){\line(1,0){2}}
\put(6,33){\line(1,0){1}}\put(11,33){\line(1,0){1}}\put(15,31){\line(1,0){1}}
\put(5,23.5){\line(1,0){1}}\put(10,24.5){\line(1,0){1}}\put(15,24.5){\line(1,0){1}}
\put(19,22.5){\line(1,0){1}}\put(5,8){\line(1,0){1}}\put(12,10){\line(1,0){1}}
\put(-0.15,0.85){$\bullet$}\put(-0.15,2.85){$\bullet$}\put(0.85,2.85){$\bullet$}\put(1.85,2.85){$\bullet$}
\put(-0.15,6.85){$\bullet$}\put(-0.15,7.85){$\bullet$}\put(-0.15,8.85){$\bullet$}
\put(0.85,8.85){$\bullet$}\put(1.85,8.85){$\bullet$}
\put(-0.15,14.35){$\bullet$}\put(-0.15,16.35){$\bullet$}\put(0.85,16.35){$\bullet$}\put(1.85,16.35){$\bullet$}
\put(-0.15,21.35){$\bullet$}\put(-0.15,22.35){$\bullet$}\put(-0.15,23.35){$\bullet$}
\put(0.85,23.35){$\bullet$}\put(1.85,23.35){$\bullet$}
\put(-0.15,29.85){$\bullet$}\put(-0.15,30.85){$\bullet$}\put(-0.15,31.85){$\bullet$}
\put(0.85,31.85){$\bullet$}\put(1.85,31.85){$\bullet$}
\put(4.85,30.85){$\circ$}\put(5.85,31.85){$\bullet$}\put(5.85,32.85){$\bullet$}
\put(9.85,30.85){$\circ$}\put(10.85,31.85){$\bullet$}\put(11.85,32.85){$\bullet$}
\put(15.85,31.85){$\bullet$}\put(16.85,32.85){$\bullet$}
\put(3.85,22.35){$\circ$}\put(4.85,24.35){$\bullet$}
\put(8.85,22.35){$\circ$}\put(9.85,23.35){$\circ$}\put(9.85,24.35){$\bullet$}
\put(13.85,22.35){$\circ$}\put(14.85,23.35){$\circ$}\put(15.85,24.35){$\bullet$}
\put(19.85,23.35){$\circ$}\put(20.85,24.35){$\bullet$}
\put(3.85,15.35){$\bullet$}\put(3.85,16.35){$\bullet$}\put(4.85,17.35){$\bullet$}
\put(5.85,8.85){$\circ$}\put(6.85,9.85){$\circ$}
\put(10.85,7.85){$\circ$}\put(11.85,8.85){$\circ$}
\put(5.85,3.85){$\circ$}\put(4.85,1.85){$\bullet$}\put(4.85,2.85){$\bullet$}
\put(0.1,29.2){$\gam,\
S$}\put(5,29.2){$W=2$}\put(10,29.2){$W=1$}\put(15,29.2){$W=1$}
\put(0.3,20.7){$\gam,\ S
$}\put(4,20.7){$W=1$}\put(9,20.7){$W=1$}\put(14,20.7){$W=1$}\put(19,20.7){$W=1$}
\put(0.3,13.7){$\gam,\ S$}\put(3,13.7){$W=-2$}
\put(0.3,6.2){$\gam, \ S$}\put(5,6.2){$W=1$}\put(11,6.2){$W=1$}
\put(0.3,0.2){$\gam,\ S$}\put(4,0.2){$W=-2$} \put(1,28.2){${\cal
R}=\{r'=2+3,\ r''=0+1+0\},\quad\widetilde m=4,\
n=3$}\put(3,27){$\text{(a) Case}\quad r'=5,\ r''=1,\ W=4$}
\put(2,19.7){${\cal R}=\{r'=1+2,\ r''=0+2+0\},\quad\widetilde
m=4,\ n=3$}\put(7,16.5){${\cal R}=\{r'=1+2,\
r''=1+0+1\}$}\put(7,15.5){$\widetilde m=3,\
n=3$}\put(3,12.5){$\text{(b) Case}\quad r'=3,\ r''=2,\
W=2$}\put(15,9){${\cal R}=\{r'=0+1,\
r''=0+3+0\}$}\put(15,8){$\widetilde m=4,\ n=3$}\put(9,3){${\cal
R}=\{r'=0+1,\ r''=1+1+1\}$}\put(9,2){$m=3,\
n=3$}\put(3,-1){$\text{(c) Case}\quad r'=1,\ r''=3,\ W=0$}
\end{picture}
\caption{Admissible paths, graphs, and markings,
III}\label{figsp10}
\end{figure}

\subsubsection{The case of
$\Sig=\bbS^2_{0,2}$}\label{secspn4} Let ${\cal R}$ be a splitting
of $r'$ and $r''$ as in (\ref{espn3}) so that
\begin{equation}r'_2+r''_{2,1}+2r''_{2,2}=m:=
|(\partial\Del)_+|+n_0\label{espn22}\end{equation} with some $0\le
n_0\le d_1$. We remark that, for $r''=0$, necessarily $n_0=0$ (see
Remark \ref{rsp1} in section \ref{secspn8} below), and the
splitting ${\cal R}$ turns into (\ref{espn53}).

An ${\cal R}$-admissible lattice path in $\Del$ is a map
$\gam:[0,m]\to\Del$ such that
\begin{enumerate}
\item[($\gam${\bf1'})] image of $\gam$ lies in ${\cal
B}_+$;\item[($\gam${\bf2'})] $\gam(0)$ and $\gam(m)$ are the two
endpoints of $(\partial\Del)_+$;\item[($\gam${\bf3'})] the
composition of the functional $x+y$ with $\gam$ is a strongly
increasing function;\item[($\gam${\bf4'})] for any $i=1,...,m$,
the point $\gam(i)$ is integral, the function
$\gam\big|_{[i-1,i]}$ is linear, and $\gam[i-1,i]$ is a unit
length segment, whose projection on ${\cal B}$ has lattice length
$1/2$ if it is not parallel to ${\cal B}$;
\item[($\gam${\bf5'})]
$\gam([0,m])\cap{\cal B}$ is finite.
\end{enumerate} We then pick integral points $v_i$,
$i=0,1,...,\widetilde m:=r'_2+r''_{2,1}+r''_{2,2}$, on $\gam$,
satisfying conditions ($\gam${\bf1})-($\gam${\bf5}) from section
\ref{secspn1} and the condition
\begin{enumerate}\item[($\gam${\bf 6'})] each break point of $\gam$, which
belongs to a segment parallel to ${\cal B}$, is among $v_i$, $1\le
i\le \widetilde m$.\end{enumerate} Denote by $\sig_i$ the segment,
joining $v_i$ with its symmetric (with respect to ${\cal B}$)
image, $i=1,...,\widetilde m-1$.

A lattice path $\gam$ with points $v_1,...,v_{\widetilde m}$ gives
rise to a set (may be, empty) of $\gam$-admissible lattice
subdivisions $S$ of $\Del$, symmetric with respect to ${\cal B}$,
and resulting from the following construction:
\begin{enumerate}\item[({\bf S1})] the part of $\Del$ between $\gam$
and its symmetric with respect to ${\cal B}$ image $\hat\gam$ is
subdivided by the segments $\sig_i$, $i=1,...,\widetilde m-1$
(see, for instance, Figure \ref{figsp2}); \item[({\bf S2})] denote
by $v_i$, $\widetilde m<i\le\widetilde m_1$, all the break points
of $\gam$, which do not appear among $v_i$, $1\le i\le\widetilde
m$;
\item[({\bf S3})] take the first break point $v$ of $\gam$, where $\gam$ turns
in the positive direction (i.e., left), and denote by
$\sig^{(1)},\sig^{(2)}$ the minimal segments of $\gam$, emanating
from $v$ and ending up at points $v_i$, $1\le i\le\widetilde m_1$;
then take a new polygon $\Del'$ of $S$ to be either

- the parallelogram, spanned by $\sig^{(1)},\sig^{(2)}$ (see
Figure \ref{figsp12}(a)), or

- if $|\sig^{(1)}|=|\sig^{(2)}|=1$, the triangle, spanned by
$\sig^{(1)},\sig^{(2)}$, with the unit length third side parallel
to ${\cal B}$ (see Figure \ref{figsp12}(b)), or

- if $|\sig^{(1)}|=2$, $|\sig^{(2)}|=1$ (resp., $|\sig^{(1)}|=1$,
$|\sig^{(2)}|=2$), the trapeze with the unit length sides
$\sig^{(3)},\widetilde\sig$, correspondingly parallel to
$\sig^{(1)}$, ${\cal B}$ (resp., $\sig^{(2)}$, ${\cal B}$) (see
Figure \ref{figsp12}(d,e)), or

- if $|\sig^{(1)}|=|\sig^{(2)}|=2$, the pentagon with the unit
length sides $\sig^{(3)},\sig^{(4)},\widetilde\sig$ parallel to
$\sig^{(1)},\sig^{(2)},{\cal B}$, respectively (see Figure
\ref{figsp12}(f));
\item[({\bf S4})] if $\Del'\not\subset\Del$ we stop the procedure,
if $\Del'\subset\Del$, we introduce the new lattice path $\gam_1$,
replacing $\sig^{(1)},\sig^{(2)}$ in $\gam$ by
$(\partial\Del')_+$, and denote the new break points of $\gam_1$
by $v_i$, $\widetilde m_1<i\le\widetilde m_2$; \item[({\bf S5})]
having a lattice path $\gam_l$ and the points $v_i$, $1\le
i\le\widetilde m_{l+1}$, we perform the above steps ({\bf S3}),
({\bf S4}), and proceed inductively until the construction stops;
\item[({\bf S6})] in case $\gam_l=(\partial\Del)_+$, we reflect
the obtained polygons with respect to ${\cal B}$ and obtain a
$\gam$-admissible subdivision $S$ of $\Del$.
\end{enumerate} Observe that such a subdivision $S$ of $\Del$ is
dual to a plane tropical curve $A$ (for the definition see, for
example, \cite{M3}, section 3.4, or section \ref{secspn6} in the
present paper below). The tropical curve $A$ has exactly
$\widetilde m$ vertices on ${\cal B}$, which we denote by
$w_1,...,w_{\widetilde m}$ in the growing coordinate order. These
vertices are dual to the polygons of $S$, symmetric with respect
to ${\cal B}$. There are $n_0$ polygons $\widetilde\Del_i$, $1\le
i\le n_0$, of $S$ in ${\cal B}_+$, different from parallelograms,
and the dual vertices of $A$ in ${\cal B}_+$ we denote by $w^+_i$,
$1\le i\le n_0$. Each of them is joined in $A$ with two of the
points $w_1,...,w_{\widetilde m}$ by segments.

Then we construct a $(\gam,{\cal R},S)$-admissible graph $G$,
starting with an auxiliary subgraph $G'$, whose connected
components are segments or points $G'_j=[(a_j,j),(b_j,j)]$ lying
on the lines $y=j$ with $j=1,...,n:=r'_1+2r''_1+r''_{2,1}$,
respectively, equipped with positive integer weights $w(G'_j)$,
satisfying conditions ({\bf G1})-({\bf G2}), ({\bf G4})-({\bf G6})
and the conditions
\begin{enumerate}\item[({\bf G7'})] for any $i=1,...,\widetilde m-1$,
relation (\ref{esp21}) holds true.
\end{enumerate}

Next we introduce some more vertices of $G$, taking one vertex
$\varphi_i$ for all $i=1,...,\widetilde m$ such that $|\gam_i|=1$,
and taking two vertices $\varphi_{i,1},\varphi_{i,2}$ for all
$i=1,...,\widetilde m$ such that $|\gam_i|=2$, and define new arcs
of $G$ following rules ({\bf G8}), ({\bf G9}) and keeping
condition ({\bf G10}) of section \ref{secspn1}.

Another set of vertices of $G$ is produced when assigning to each
vertex $w^+_i$ a pair of vertices
$\varphi^+_{i,1},\varphi^+_{i,2}$ of $G$. Then we introduce
additional arcs:
\begin{enumerate}\item[({\bf G12})] if $|\gam_i|=1$ and
the vertex $w_i$ of $A$ is joined by a segment with a vertex
$w^+_j$, then we connect the vertex $\varphi_i$ of $G$ both with
$\varphi^+_{j,1}$ and $\varphi^+_{j,2}$ by arcs;
\item[({\bf G13})] if $|\gam_i|=2$ and the vertex $w_i$ of $A$ is joined
by a segment with only one vertex $w^+_j$, $1\le j\le n_0$, then
either we connect by arcs the vertex $\varphi_{i,1}$ of $G$ with
$\varphi^+_{j,1}$, and the vertex $\varphi_{i,2}$ with
$\varphi^+_{j,2}$, or we connect $\varphi_{i,1}$ with
$\varphi^+_{j,2}$, and $\varphi_{i,2}$ with $\varphi^+_{j,1}$;
\item[({\bf G14})] if $|\gam_i|=2$ and the vertex $w_i$ of $A$ is joined by segments with
two vertices $w^+j$, $w^+_k$, then either we connect the vertex
$\varphi_{i,1}$ both with $\varphi^+_{j,1}$ and $\varphi^+_{k,2}$,
respectively, the vertex $\varphi_{i,2}$ both with
$\varphi^+_{j,2}$ and $\varphi^+_{k,1}$, or we connect the vertex
$\varphi_{i,1}$ both with $\varphi^+_{j,2}$ and $\varphi^+_{k,1}$,
respectively, the vertex $\varphi_{i,2}$ both with
$\varphi^+_{j,1}$ and $\varphi^+_{k,2}$.
\end{enumerate} Our final requirement to $G$ is that is must be a tree
(condition ({\bf G11}) of section \ref{secspn1}).

A marking of a $(\gam,{\cal R},S)$-admissible graph $G$ is a pair
of integer vectors \mbox{$\overline s'=(s'_1,...,s'_{r'_1})$} and
$\overline s''=(s''_1,...,s''_{r''_1})$, satisfying conditions
({\bf s1})-({\bf s4}) of section \ref{secspn1}.

We define the Welschinger number as $W({\cal R},\gam,G,\overline
s',\overline s'')=0$ if at least one weight $w(G'_i)$,
$2r''_1<i\le n$, is even, and otherwise as $$W({\cal
R},\gam,S,G,\overline s',\overline s'')=(-1)^a\cdot
2^{b+c}\cdot\prod_{j=1}^{r''_1}(w(G'_{2j}))^2\cdot\prod_{j=2r''_1+1}^{2r''_1+r''_{2,1}}
w(G'_j)$$ \begin{equation}\times\prod_{i=0}^{\widetilde m}
\left(n'_i!n''_i!\alp_i^{-1}\bet_i^{-1}\right)\cdot\prod_{k=1}^{n_0}\Area(\widetilde\Del_k)\
,\label{espn62}\end{equation}
$$\alp_i=\prod_{\renewcommand{\arraystretch}{0.6}
\begin{array}{c}
\scriptstyle{0\le d\le e\le\widetilde m}\\
\scriptstyle{f=1,3,5,...}
\end{array}}
n'_{i,d,e,f}!,\quad\bet_i=\prod_{\renewcommand{\arraystretch}{0.6}
\begin{array}{c}
\scriptstyle{0\le d\le e\le\widetilde m}\\
\scriptstyle{f=1,2,3,...}
\end{array}}
n''_{i,d,e,f}!\ ,$$ where \begin{itemize}\item $a$, $b$ are as in
formula (\ref{espn61}); \item $c$ is the number of the polygons
$\Del'$ in the subdivision $S$, which are symmetric with respect
to ${\cal B}$ and whose boundary part $(\partial\Del')_+$ is not a
segment;
\item $n'_i=\#\{j\ :\ s'_j=i,\ 1\le j\le r'_1\}$, $n''_i=\#\{j\ :\
s''_j=i,\ 1\le j\le r''_1\}$;\item $n'_{i,d,e,f}=\#\{j\ :\
s'_j=i,\ a_j=d,\ b_j=e,\ w(G'_j)=f\}$;
\item $n''_{i,d,e,f}=\#\{j\ :\ s''_j=i,\ a_{2j}=d,\
b_{2j}=e,\ w(G'_{2j})=f\}$; \item $\Area(\widetilde\Del_k)$ is the
following: $\widetilde\Del_k$ uniquely splits into the Minkowski
sum of a triangle with $0$, $1$, or $2$ segments, and we put
$\Area(\widetilde\Del_k)$ to be the lattice area of that triangle.
\end{itemize}

\begin{theorem}\label{tsp4}
In the notation of section \ref{secspn4}, if $\Sig=\bbS^2_{0,2}$,
and the positive integers $r',r''$ satisfy (\ref{espn1}), then
\begin{equation}W_{r''}(\Sig,{\cal L}_\Del)=\sum W({\cal R},\gam,S,G,\overline s',\overline s'')\
,\label{espn40}\end{equation} where the sum ranges over splittings
${\cal R}$ of $r',r''$, satisfying (\ref{espn2}), all ${\cal
R}$-admissible lattice paths $\gam$, all $\gam$-admissible
subdivisions $S$ of $\Del$, all $(\gam,{\cal R},S)$-admissible
graphs $G$, and all markings $\overline s',\overline s''$ of $G$,
subject to conditions, specified in section \ref{secspn4}.
\end{theorem}

As example we consider the linear system $|{\cal L}_\Del|$ with
the hexagon $\Del$ shown in Figure \ref{fn2}(d) for $d=2$,
$d_1=1$. Here $r'+2r''=5$, and the integration with respect to
Euler characteristic (cf. section \ref{secspn2}, case (A)) gives
$W_{r''}(\bbS^2_{0,2},{\cal L}_\Del)=4-2r''$, $r''=1,2$. In Figure
\ref{figsp11} we demonstrate how to obtain these answers from
formula (\ref{espn40}), listing admissible paths $\gam$,
subdivisions $S$ of $\Del$, and graphs $G'$ (here the marked
points are denoted by bullets and the non-marked one-point
components of $G'$ are denoted by circles).

\begin{figure}
\setlength{\unitlength}{0.6cm}
\begin{picture}(19,23)(-3,0)
\thinlines \put(0,2){\vector(1,0){3}}\put(4.5,2){\vector(1,0){2}}
\put(0,2){\vector(0,1){3}}\put(4.5,2){\vector(0,1){3}}
\put(0,7.5){\vector(1,0){3}}\put(4.5,7.5){\vector(1,0){2}}
\put(0,7.5){\vector(0,1){3}}\put(4.5,7.5){\vector(0,1){3}}
\put(0,14){\vector(1,0){3}}\put(4.5,14){\vector(1,0){2}}
\put(0,14){\vector(0,1){3}}\put(4.5,14){\vector(0,1){3}}
\put(0,19.5){\vector(1,0){3}}\put(4.5,19.5){\vector(1,0){2}}
\put(0,19.5){\vector(0,1){3}}\put(4.5,19.5){\vector(0,1){3}}
\put(2,3){\line(0,1){1}}\put(0,15){\line(1,1){1}}
\put(1,2){\line(0,1){1}}\put(1,3){\line(1,0){1}}
\put(0,3){\line(1,1){1}}\put(1,2){\line(1,1){1}}
\put(1,7.5){\line(1,1){1}}\put(1,7.5){\line(-1,1){1}}
\put(2,8.5){\line(0,1){1}}\put(2,8.5){\line(-1,1){1}}
\put(1,14){\line(0,1){1}}\put(1,14){\line(1,1){1}}
\put(1,15){\line(1,0){1}}\put(2,15){\line(0,1){1}}
\put(2,15){\line(-1,1){1}}\put(1,19.5){\line(-1,1){1}}
\put(1,19.5){\line(1,1){1}}\put(2,20.5){\line(-1,1){1}}
\put(2,20.5){\line(0,1){1}} \thicklines \put(0,2){\line(0,1){1}}
\put(0,3){\line(1,0){1}}\put(1,3){\line(0,1){1}}\put(1,4){\line(1,0){1}}
\put(0,7.5){\line(0,1){1}}\put(0,8.5){\line(1,1){1}}\put(1,9.5){\line(1,0){1}}
\put(0,14){\line(0,1){1}}\put(0,15){\line(1,0){1}}\put(1,15){\line(0,1){1}}
\put(1,16){\line(1,0){1}}\put(0,19.5){\line(0,1){1}}\put(0,20.5){\line(1,1){1}}
\put(1,21.5){\line(1,0){1}}
\put(5.35,21.35){$\bullet$}\put(5.35,14.85){$\bullet$}
\put(5.35,20.35){$\circ$}\put(5.35,8.35){$\circ$}\put(5.35,9.35){$\circ$}
\put(5.35,2.85){$\circ$}\put(0.5,1){$\gam,\
S$}\put(0.5,6.5){$\gam,\ S$}\put(0.5,13){$\gam,\
S$}\put(0.5,18.5){$\gam,\ S$}\put(4.6,1){$G',\ \overline
s$}\put(4.6,6.5){$G',\ \overline s$}\put(4.6,13){$G',\ \overline
s$}\put(4.6,18.5){$G',\ \overline
s$}\put(8,3){$\begin{cases}&{\cal R}=\{r'=0+1,\ r''=0+1+1\}\\
&\widetilde m=1,\quad n=1\\
&W=-2\end{cases}$}\put(8,8.5){$\begin{cases}&{\cal R}=\{r'=0+1,\ r''=0+2+0\}\\
&\widetilde m=1,\quad n=2\\
&W=2\end{cases}$}\put(8,15){$\begin{cases}&{\cal R}=\{r'=1+2,\ r''=0+0+1\}\\
&\widetilde m=1,\quad n=1\\
&W=-2\end{cases}$}\put(8,20.5){$\begin{cases}&{\cal R}=\{r'=1+2,\ r''=0+1+0\}\\
&\widetilde m=1,\quad n=2\\
&W=4\end{cases}$}\put(1,0){$\text{(b) Case}\ r'=1,\ r''=2,\
W=0$}\put(1,12){$\text{(a) Case}\ r'=3,\ r''=1,\ W=2$}
\end{picture}
\caption{Admissible paths, graphs, and markings,
IV}\label{figsp11}
\end{figure}

\section{Tropical limits of real rational curves on non-standard real toric Del Pezzo
surfaces}\label{secspn7}

\subsection{Preliminaries}\label{secspn6} Here we recall definitions and a few facts about
tropical curves and tropical limits of algebraic curves over a
non-Archimedean field, presented in \cite{I,K,M2,M3,RST,ShP,ShW}
in more details.

By Kapranov's theorem the amoeba $A_C$ of a curve $C\in|{\cal
L}_\Del|_\K$, given by an equation
\begin{equation}f(x,y):=\sum_{(i,j)\in\Del}a_{ij}x^iy^j=0,\quad a_{ij}\in\K,\ (i,j)\in\Del\cap\Z^2,
\label{enn25}\end{equation} with the Newton polygon $\Del$, is the
corner locus of the convex piece-wise linear function
\begin{equation}N_f(x,y)=\max_{(i,j)\in
\Del\cap\Z^2}(xi+yj+\val(a_{ij})),\quad x,y\in\R\
.\label{enn1}\end{equation} In particular, $A_C$ is a planar graph
with all vertices of valency $\ge 3$.

Take the convex polyhedron
$$\widetilde\Del=\{(i,j,\gam)\in\R^3\ :\ \gam\ge-\val(a_{ij}),\quad (i,j)\in
\Del\cap\Z^2\}$$ and define the function
\begin{equation}\nu_f:\Del\to\R,\quad\nu_f(x,y)=\min\{\gam\ :\
(x,y,\gam)\in\widetilde\Del\}\ .\label{enn4}\end{equation} This is
a convex piece-wise linear function, whose linearity domains form
a subdivision $S_C$ of $\Del$ into convex lattice polygons
$\Del_1,...,\Del_N$. The function $\nu_f$ is Legendre dual to
$N_f$, and the subdivision $S_C$ is combinatorially dual to the
pair $(\R^2,A_C)$. Clearly, $A_C$ and $S_C$ do not depend on the
choice of a polynomial $f$ defining the curve $C$.

We define the {\it tropical curve}, corresponding to the algebraic
curve $C$, as the weighted graph, supported at $A_C$, i.e., the
non-Archimedean amoeba $A_C$, whose edges are assigned the weights
equal to the lattice lengths of the dual edges of $S_C$. The
subdivision $S_C$ can be uniquely restored from the tropical curve
$A_C$.

By the {\it tropical limit} of a curve $C$ given by (\ref{enn25})
we call a pair $(A_C,\{C_1,...,C_N\})$, where $C_k$, $1\le k\le
N$, is a complex curve on the toric surface $\Tor(\Del_k)$,
associated with a polygon $\Del_k$ from the subdivision $S_C$, and
is defined by an equation
$$f_k(x,y):=\sum_{(i,j)\in\Del_k}a^0_{ij}x^iy^j=0\ ,$$ where
$a_{ij}(t)=(a^0_{ij}+O(t^{>0}))\cdot t^{\nu_f(i,j)}$ is the
coefficient of $x^iy^j$ in $f(x,y)$. We call $C_1,...,C_N$ {\it
limit curves}. Their geometrical meaning is as follows (cf.
\cite{ShP}, section 2). By a parameter change $t\mapsto t^M$,
$M>>1$, we can make all the exponents of $t$ in the coefficients
$a_{ij}=a_{ij}(t)$ of $f$ integral, and make the function $\nu_f$
integral-valued at integral points. The toric threefold
$Y=\Tor(\widetilde\Del)$ fibers over $\C$ so that $Y_t$, $t\ne 0$,
is isomorphic to $\Tor(\Del)$, and $Y_0$ is the union of
$\Tor(\Del_k)$ attached to each other as the polygons of the
subdivision $S_C$. Equation (\ref{enn25}) defines an analytic
surface $C$ in $Y$ such that the curves $C^{(t)}=C\cap Y_t$,
$0<|t|<\eps$, form an equisingular family, and $C^{(0)}=C\cap
Y_0=C_1\cup...\cup C_N$, where $C_k=C\cap\Tor(\Del_k)$.

Now let $\Sig=\bbS^2$, $\bbS^2_{1,0}$, $\bbS^2_{2,0}$, or
$\bbS^2_{0,2}$ be defined over the field $\K$, and let $\Del$ be
one of the respective polygons shown in Figure \ref{fn2}(a-d). Fix
some non-negative integers $r',r''$ satisfying (\ref{enn36}) and
pick a configuration $\op$ of $-c_1({\cal L}_\Del)K_{\Sig}-1$
distinct points in $(\K^*)^2\subset\Sig$ such that
$\op=\op'\cup\op''$ with
\begin{equation}\begin{cases}\op'=\{\bp'_1,...,\bp'_{r'}\}\subset\Sig(\K_\R),\\
\op''=\{\bp''_{1,1},\bp''_{1,2},...,\bp''_{r'',1},\bp''_{r'',2}\}
\subset\Sig(\K)\backslash\Sig(\K_R),\\
\conj(\bp''_{i,1})=\bp''_{i,2},\quad
i=1,...,r''.\end{cases}\label{espn21}\end{equation} Since the
anti-holomorphic involution acts on $(\K^*)^2\subset\Sig(\K)$ by
$\conj(\xi,\eta)=(\overline\eta,\overline\xi)$, we have
$\bp'_i=(\xi_i(t),\overline\xi_i(t))$, $i=1,...,r'$, and
$\bp''_{i,1}=(\eta_i(t),\zeta_i(t))$,
$\bp''_{i,2}=(\overline\zeta_i(t),\overline\eta_i(t))$,
$i=1,...,r''$. In particular, the configuration $\ox'=\val(\op')$
lies on ${\cal B}$, and the configuration $\ox''=\val(\op'')$ is
symmetric with respect to ${\cal B}$. We assume $\op$ to be
generic in $\Omega_{r''}(\Sig(\K),{\cal L}_\Del)$ and such that
the configuration $\ox=\ox'\cup\ox''$ consists of $r'+r''$ generic
distinct points on ${\cal B}$: $$\ox'=\{\bx'_i\ :\
\bx'_i=\val(\bp'_i),\ i=1,...,r'\}\ ,$$ $$\ox''=\{\bx''_i\ :\
\bx''_i=\val(\op''_{i,1})=\val(\bp''_{i,2}),\ i=1,...,r''\}\ .$$
In addition, we require the following property:

\begin{enumerate}\item[($\bx${\bf1})] the points of $\ox$ are ordered on
${\cal B}$ as
$$\bx''_1\prec...\prec\bx''_{r''}\prec\bx'_1\prec...\prec\bx'_{r'}$$
by the growing sum of coordinates, and, moreover, the distance
between any pair of neighboring points is much larger than that
for the preceding pair.
\end{enumerate}

Let $C\in|{\cal L}_\Del|_\K$ be a real rational curve, passing
through $\op$. We can define $C$ by an equation (\ref{enn25}) with
$a_{ji}=\overline a_{ij}$, $(i,j)\in\Del$. Observe that the
tropical curve $A_C\subset\R^2$ is symmetric with respect to
${\cal B}$, and so is the dual subdivision $S_C$ of $\Del$. We
intend to describe the tropical limits $(A_C,\{C_1,...,C_N\})$ of
such curves $C$.

\subsection{Tropical limits of real rational curves on $\bbS^2$, $\bbS^2_{1,0}$,
or $\bbS^2_{2,0}$}\label{secspn3} In addition to the notation of
the preceding section, introduce the following ones:
\begin{itemize}\item let $P(S_\Del)$, $E(S_C)$, and $V(S_C)$ be the sets of the
polygons, the edges, and the vertices of $S_C$; \item $P(S_C)$
splits into the disjoint subsets $P_{\cal B}(S_C)$, containing the
polygons symmetric with respect to ${\cal B}$, and $P_+(S_C)$,
$P_-(S_C)$, consisting of the polygons contained in the
half-planes ${\cal B}_+$, ${\cal B}_-$, respectively; notice that
the polygons of $P_+(S_C)$ and $P_-(S_C)$ are in 1-to-1
correspondence by the reflection with respect to ${\cal B}$, and
that the limit curves $C_i$, $\Del_i\in P_{\cal B}(S_C)$ are real,
and the limit curves $C_j$, $\Del_j\not\in P_{\cal B}(S_C)$ are
not;
\item for a limit curve $C_i$, $1\le i\le N$, denote by $C_{ij}$,
$1\le j\le l_i$, the set of all its components (counting each one
with its multiplicity);
\item for any $\Del_k\in P(S_C)$, denote by $\Tor(\partial\Del_k)$
the union of the toric divisors $\Tor(\sig)$,
$\sig\subset\partial\Del_k$.
\end{itemize} Observe that all the curves $C_{ij}$, $1\le j\le
l_i$, $1\le i\le N$, are rational (cf. \cite{ShW}, Step 1 of the
proof of Proposition 2.1), which comes from the inequality for
geometric genera $g(C^{(t)})\ge\sum_{i,j}g(C_{ij})$ (see
\cite{DH}, Proposition 2.4, or \cite{N}. We say that $C_{ij}$ is
{\bf binomial} and write $C_{ij}\in{\cal C}_b$ if $C_{ij}$
intersects with $\Tor(\partial\Del_i)$ at precisely two points (in
this case the two toric divisors, meeting $C_{ij}$, correspond to
opposite parallel edges of $\Del_i$, and $C_{ij}$ is defined by an
irreducible binomial), otherwise we write $C_{ij}\in{\cal
C}_{nb}$.

Now we split the configuration $\ox$ as follows:
$$\ox'=\ox'_1\cup\ox'_2\cup\ox'_3,\quad\ox''=\ox''_1\cup\ox''_2\cup\ox''_3\
,$$ where \begin{itemize}\item $\ox'_1$ (resp., $\ox''_1$)
consists of the points $\bx'_i$, $1\le i\le r'$ (resp., $\bx''_i$,
$1\le i\le r''$), lying inside the edges of $A_C$ on ${\cal B}$;
\item $\ox'_2$ (resp., $\ox''_2$) consists of the points $\bx'_i$, $1\le i\le r'$
(resp., $\bx''_i$, $1\le i\le r''$), which are vertices of $A_C$;
\item $\ox'_3$ (resp., $\ox''_3$) consists of the points $\bx'_i$, $1\le i\le r'$
(resp., $\bx''_i$, $1\le i\le r''$), lying inside edges of $A_C$
orthogonal to ${\cal B}$.
\end{itemize} Furthermore, a point $\bx'_i\in\ox'_2$ (resp., $\bx''_i\in\ox''_2$) is dual to
a polygon $\Del_k\in P_{\cal B}(S_C)$. In particular, the point
$\bp'_i\in(\R^*)^2\subset\Tor(\Del_k)$ (resp., the points
$\bp''_{i,1},\bp''_{i,2}\in(\C^*)^2\subset\Tor(\Del_k)$) lies on
the limit curve $C_k$. For a point
$\bp=(\xi(t),\eta(t))\in(\K^*)^2$, we put
$\ini(\bp)=(\xi_0,\eta_0)\in(\C^*)^2$, $\xi_0,\eta_0$ being the
coefficients of the lower powers of $t$ in $\xi(t),\eta(t)$,
respectively. We then have
$\ox'_2=\ox'_{2,1}\cup\ox'_{2,2}\cup\ox'_{2,3}\cup\ox'_{2,4}$ and
$\ox''_2=\ox''_{2,1}\cup\ox''_{2,2}\cup\ox''_{2,3}\cup\ox''_{2,4}$,
where \begin{itemize}\item $\ox'_{2,1}$ (resp., $\ox''_{2,1}$)
consists of the points $\bx'_i$, $1\le i\le r'$ (resp., $\bx''_i$,
$1\le i\le r''$), such that $\ini(\bp'_i)$ lies (resp., both
$\ini(\bp''_{i,1})$ and $\ini(\bp''_{i,2})$ lie) on a real
non-binomial component $C_{kl}$ of $C_k$; we, furthermore, make
splitting $\ox''_{2,1}=\ox''_{2,1a}\cup\ox''_{2,1b}$, where a
point $\bx''_i\in\ox''_{2,1}$, dual to $\Del_k\in P_{\cal
B}(S_C)$, belongs to $\ox''_{2,1a}$ or to $\ox''_{2,1b}$ according
as $C_{kl}$ has one or at least two local branches centered along
$\Tor((\partial\Del_k)_+)$;
\item $\ox'_{2,2}$ (resp., $\ox''_{2,2}$) consists of the points
$\bx'_i$, $1\le i\le r'$ (resp., $\bx''_i$, $1\le i\le r''$), such
that the point $\ini(\bp'_i)$ lies (resp., the points
$\ini(\bp''_{i,1})$, $\ini(\bp''_{i,2})$ lie) on two distinct
conjugate components of $C_k$ which cross
$\Tor((\partial\Del_k)_+)$;
\item $\ox'_{2,3}$ (resp., $\ox''_{2,3}$) consists of the points
$\bx'_i$, $1\le i\le r'$ (resp., $\bx''_i$, $1\le i\le r''$), such
that $\ini(\bp'_i)$ lies on a real binomial component $C_{kj}$
(resp., $\ini(\bp''_{i,1})$, $\ini(\bp''_{i,2})$ lie on two
distinct conjugate binomial components $G_{kj},C_{kl}$) crossing
$\Tor((\partial\Del_k)_\perp)$.
\end{itemize}

Put $r'_i=\#(\ox'_i)$, $r''_i=\#(\ox''_i)$, $i=1,2,3$,
$r'_{2,j}=\#(\ox'_{2,j})$, $r''_{2,j}=\#(\ox''_{2,j})$,
$j=1,2,3,4$, $r''_{2,1a}=\#(\ox''_{2,1a})$,
$r''_{2,1b}=\#(\ox''_{2,1b})$.

Projecting the polygons $\Del_k\in P_{\cal B}(S_C)$ to
$(\partial\Del)_+$ in the direction orthogonal to ${\cal B}$, we
obtain
$$|\partial\Del|-|(\partial\Del)_\perp|=2|(\partial\Del)_+|=2\sum_{\renewcommand{\arraystretch}{0.6}
\begin{array}{c}
\scriptstyle{\sig\subset(\partial\Del_k)_+}\\
\scriptstyle{\Del_k\in P_{\cal B}(S_C)}
\end{array}}\pr(\sig)+4\sum_{\renewcommand{\arraystretch}{0.6}
\begin{array}{c}
\scriptstyle{\sig\in E(S_C)}\\
\scriptstyle{\sig\subset{\cal B}}
\end{array}}|\sig|$$ \begin{equation}\ge2(r'_{2,1}+2r'_{2,2}+r'_{2,3}+
r''_{2,1a}+2r''_{2,1b}+2r''_{2,2}+r''_{2,3})+2(2r'_3+4r''_3)\
,\label{espn5}\end{equation} where the equality holds only if
\begin{enumerate}\item[({\bf E1})] each polygon
$\Del_k\in P_{\cal B}(S_C)$ is dual to a point of $\ox'_2\cup\ox''_2$, and

- either $\pr((\partial\Del_k)_+)=|(\partial\Del_k)_+|=1$, the
limit curve $C_k$ may have binomial components, meeting the toric
divisors $\Tor(\sig)$, $\sig\subset(\partial\Del_k)_\perp$, and it
has one real non-multiple component, crossing the toric divisor
$\Tor((\partial\Del_k)_+)$,

- or $\pr((\partial\Del_k)_+)=|(\partial\Del_k)_+|=2$, the limit
curve $C_k$ may have binomial components, meeting the toric
divisors $\Tor(\sig)$, $\sig\subset(\partial\Del_k)_\perp$, and it
has two distinct non-multiple conjugate components, meeting the
toric divisors $\Tor((\partial\Del_k)_+)$,

- or $\pr((\partial\Del_k)_+)=|(\partial\Del_k)_+|=2$, the limit
curve $C_k$ may have binomial components, meeting the toric
divisors $\Tor(\sig)$, $\sig\subset(\partial\Del_k)_\perp$, and it
has one real non-multiple component, crossing the toric divisors
$\Tor((\partial\Del_k)_+)$;
\item[({\bf E2})] any point $\ini(\bp'_i)$ (resp., $\ini(\bp''_{i,1})$ or $\ini(\bp''_{i,2})$)
with $\bx'_i\in\ox'_2$ (resp., $\bx''_i\in\ox''_2$), dual to a
polygon $\Del_k\in P_{\cal B}(S_C)$, lies precisely on one
component of the limit curve $C_k$.
\end{enumerate}

Next, inequality (2.7) from \cite{ShW}, in our situation, reads
$$\sum_{\Del_k\in P_{\cal
B}(S_C)}\sum_{j=1}^{l_k}(B(C_{kj})-2)+\sum_{\Del_k\in
P(S_C)\backslash P_{\cal B}(S_C)}\sum_{j=1}^{l_k}(B(C_{kj})-2)$$
\begin{equation}\le B(\partial\Del)-2\le|\partial\Del|-2\
,\label{espn6}\end{equation} where $B(C_{kj})$ is the number of
the local branches of the curve $C_{kj}$ centered along
$\Tor(\Del_k)$, and $B(\partial\Del)$ is the number of the local
branches of all the curves $C_k$ (counting multiplicities)
centered along the toric divisors $\Tor(\sig)$,
$\sig\subset\Del_k\cap\partial\Del$, $k=1,...,N$. The equalities
in (\ref{espn6}) hold only if \begin{enumerate}\item[({\bf E3})]
local branches of the curves $C_k$ centered along $\Tor(\sig)$,
$\sig\subset\Del_k\cap\partial\Del$, $1\le k\le N$, are
nonsingular and transverse to $\Tor(\sig)$; \item[({\bf E4})] for
any edge $\sig=\Del_k\cap\Del_l$ and any point $p\in C_k\cap
C_l\subset\Tor(\sig)$, the numbers of local branches of
$(C_k)_{\red}$ and $C_l)_{\red}$ at $p$ coincide (cf. \cite{ShP},
Remark 3.4); \item[({\bf E5})] no singular point of any curve
$(C_k)_{\red}$ in $(\C^*)^2\subset\Tor(\Del_k)$, $k=1,...,N$, is
smoothed up in the deformation $C^{(t)}$, $t\in(\C,0)$.
\end{enumerate}

To proceed with estimations, we introduce some auxiliary objects.
For any edge $\sig\subset\bigcup_{\Del_k\in P_{\cal
B}(S_C)}(\partial\Del_k)_\perp$, there is a canonical
identification $\pi_\sig:\PP^1\to\Tor(\sig)$. For such an edge
$\sig$, put $\Phi_\sig=\Tor(\sig)\cap C^{(0)}$ and define
$\Phi=\bigcup_{\sig}\pi_\sig^{-1}(\Phi_\sig)\subset\C^*$. Denote
by $\widetilde\Phi$ the set of the local branches of the
non-binomial components of the curves $(C_k)_{\red}$, $\Del_k\in
P_{\cal B}(S_C)$, centered at $\bigcup_\sig\Phi_\sig$. We claim
that
$$\#(\widetilde\Phi)\ge2(r'_1+2r''_1+r'_{2,3}+2r''_{2,3}+r''_{2,1a})-B_0$$ \begin{equation}\ge
2(r'_1+2r''_1+r'_{2,3}+2r''_{2,3}+r''_{2,1a})-|(\partial\Del)_\perp|\
,\label{espn8}\end{equation} where $$B_0=
\#\left(C^{(0)}\cap\Tor\left(\bigcup_{\Del_k\in P_{\cal
B}(S_C)}(\partial\Del_k)_\perp\cap\partial\Del\right)\right)$$
Indeed, observe that
\begin{itemize}\item each point $\bp\in\op$ such that
$\val(\bp)\in\ox'_1\cup\ox''_1\cup\ox'_{2,3}\cup\ox''_{2,3}$
defines a point in $\Phi$, and by our choice all these points are
distinct and generic; \item binomial components of the curves
$C_k$, $\Del_k\in P_{\cal B}(S_C)$, may join only the points in
$\bigcup_\sig\Phi_\sig$ with the same image in $\Phi$.
\end{itemize} Furthermore, let $I\subset\{1,2,...,N\}$ consist of all $k$ such that
the polygon $\Del_k$ belongs to $P_{\cal B}(S_C)$ and is dual to a
point of $\ox''_{2,1a}$. Let $C^{nb}_k$ be the non-binomial
component of the limit curve $C_k$, where $k\in I$, and let
$\Del_k^{nb}$ be its Newton polygon. Since
$|(\partial\Del_k^{nb})_+|=|(\partial\Del_k)_+|=1$, the space of
rational curves in the linear system $|{\cal L}_{\Del_k^{nb}}|$ on
the surface $\Tor(\Del_k)$, passing through the points
$\ini(\bp''_{i,1}),\ini(\bp''_{i,2})$, where
$\bx''_i\in\ox''_{2,1a}$ is dual to $\Del_k$, is of dimension
$|(\partial\Del_k^{nb})_\perp|-1$. That is, if $s$ is the total
number of the local branches of all the curves $C_k^{nb}$,
centered along $\Tor((\partial\Del_k)_\perp)$, $k\in I$, no more
than $s-\#(I)$ of them can be chosen in a generic position, and
hence the bound (\ref{espn8}) follows. Next, we derive
$$\sum_{\Del_k\in P_{\cal
B}(S_C)}\sum_{j=1}^{l_k}(B(C_{kj})-2)\ge\#(\widetilde\Phi)+2r''_{2,1b}$$
\begin{equation}\ge
2(r'_1+2r''_1+r'_{2,3}+2r''_{2,3}+r''_{2,1a}+r''_{2,1b})-|(\partial\Del)_\perp|\
, \label{espn9}\end{equation} where the equality holds only if
\begin{enumerate}\item[({\bf E6})] the sides of $\Del$, orthogonal to
${\cal B}$, are sides of some polygons $\Del_k\in P_{\cal
B}(S_C)$, and all the intersection points of $C^{(0)}$ with
$\Tor((\partial\Del)_\perp)$ are non-singular and transversal;
\item[({\bf E7})] $\#(\Phi)=2(r'_1+2r''_1+r'_{2,3}+2r''_{2,3}+r''_{2,1a})$, all
the points of $\bigcup_\sig\Phi_\sig$ with the same image in
$\Phi$ are joined by binomial components of the curves $C_k$,
$\Del_k\in P_{\cal B}(S_C)$, into one connected component, and any
curve $(C_k)_{\red}$ with $\Del_k\in P_{\cal B}(S_C)$ is unibranch
at each point of $\Phi_\sig$, $\sig\subset(\partial\Del_k)_\perp$.
\end{enumerate}

Now, from (\ref{espn1}) we derive
\begin{equation}2|\partial\Del|-2=2r'+4r''=2(r'_1+r'_{2,1}+r'_{2,2}+r'_{2,3}+r'_3)$$
$$+4(r''_1+r''_{2,1a}+r''_{2,1b}+r''_{2,2}+r''_{2,3}+r''_3)\
,\label{espn16}\end{equation} which after subtracting inequality
(\ref{espn5}) gives
\begin{equation}|\partial\Del|+|(\partial\Del)_\perp|-2\le2r'_1+4r''_1+2r''_{2,1a}+2r''_{2,3}-2r'_{2,2}-
2r'_3-4r''_3\ .\label{espn10}\end{equation} On the other hand,
(\ref{espn6}) and (\ref{espn9}) yield
$$|\partial\Del|+B_0-2\ge\sum_{\Del_k\in
P(S_C)\backslash P_{\cal B}(S_C)}\sum_{j=1}^{l_k}(B(C_{kj})-2)$$
$$+2r'_1+4r''_1+2r'_{2,3}+2r''_{2,1a}+2r''_{2,1b}+4r''_{2,3}\ ,$$
which together with (\ref{espn10}) and
$B_0\le|(\partial\Del)_\perp|$ results in
\begin{equation}B(C_{kj})=2\quad\text{for all}\quad\Del_k\in P(S_C)\backslash P_{\cal B}(S_C),
\ j=1,...,l_k\ ,\label{espn11}\end{equation}
\begin{equation}r'_{2,3}=r''_{2,3}=r''_{2,1b}=r'_{2,2}=r'_3=r''_3=0\
,\label{espn12}\end{equation} and implies the equalities in
(\ref{espn5}), (\ref{espn6}), (\ref{espn9}), that is all the
conditions ({\bf E1})-({\bf E7}) hold true as well as
\begin{enumerate}\item[({\bf E8})] for any $\Del_k\in P(S_C)\backslash
P_{\cal B}(S_C)$, the curve $C_k$ consists of only binomial
components (see (\ref{espn11})).\end{enumerate}

\begin{remark}\label{rspn1} From conditions ({\bf E1})-({\bf E9}) and
equalities (\ref{espn11}), (\ref{espn12}), one can easily derive
that $\bigcup_{\Del_k\in P_{\cal B}(S_C)}(\partial\Del_k)_+$ with
the respective vertices of the polygons $\Del_k\in P_{\cal
B}(S_C)$ forms a connected lattice path $\gam$ in $\Del$,
satisfying conditions ($\gam${\bf1})-($\gam${\bf5}). Furthermore,
due to ({\bf E5}) and ({\bf E8}), $\gam$ has no intersections with
${\cal B}$ besides its endpoints (since, otherwise, the curves
$C{(t)}$, $t\ne 0$, would be reducible). At last, placing the
configuration $\ox$ on ${\cal B}$ so that $\ox''$ will precede
$\ox'$, we get ($\gam${\bf4}).
\end{remark}

\subsection{Tropical limits of real rational curves on
$\Sig=\bbS^2_{0,2}$}\label{secspn8} In the considered situation,
the argument of the preceding section leads to inequalities
(\ref{espn6}), (\ref{espn8}), and (\ref{espn9}), the latter one
turning, due to \mbox{$(\partial\Del)_\perp=\emptyset$}, $B_0=0$
(see Figure \ref{fn2}(d)), into \begin{equation}\sum_{\Del_k\in
P_{\cal B}(S_C)}\sum_{j=1}^{l_k}(B(C_{kj})-2)\ge
2(r'_1+2r''_1+r'_{2,3}+2r''_{2,3}+r''_{2,1a}+r''_{2,1b})\ ,
\label{espn14}\end{equation} Inequality (\ref{espn5}) will be
replaced by another relations for the parameter $B_+$, equal to
the total number of the local branches of the curves $(C_k)_\red$,
$\Del_k\in P_{\cal B}(S_C)$, centered on the divisors
$\Tor((\Del_k)_+)$, and the local branches of the curves
$(C_j)_\red$, $\Del_j\in P_+(S_C)$, centered on the divisors
$\Tor(\sig)$, $\sig\subset{\cal B}$.

We have two possibilities.

\medskip

{\bf Case 1}. Assume that
\begin{equation}B_+\le|(\partial\Del)_+|=2d-d_1\
.\label{espn42}\end{equation} In the notation of section
\ref{secspn3}, this yields
\begin{equation}|\partial\Del|=2|(\partial\Del)_+|\ge2(r'_{2,1}+2r'_{2,2}+r'_{2,3}+
r''_{2,1a}+2r''_{2,1b}+2r''_{2,2}+r''_{2,3}+r'_3+2r''_3)\label{espn15}\end{equation}
with an equality only if \begin{enumerate}\item[({\bf E1'})] for
each polygon $\Del_k\in P_{\cal B}(S_C)$,

- either $|(\partial\Del_k)_+|=1$, $\Del_k$ is dual to a point
from $\ox'_{2,1}\cup\ox''_{2,1a}$, the limit curve $C_k$ may have
binomial components, meeting the toric divisors $\Tor(\sig)$,
$\sig\subset(\partial\Del_k)_\perp$, and it has one real (possibly
multiple) component, crossing the toric divisor
$\Tor((\partial\Del_k)_+)$,

- or $|(\partial\Del_k)_+|=2$, $\Del_k$ is dual to a point from
$\ox''_{2,2}$, the limit curve $C_k$ may have binomial components,
meeting the toric divisors $\Tor(\sig)$,
$\sig\subset(\partial\Del_k)_\perp$, and it has two distinct
(possibly multiple) conjugate components, meeting the toric
divisors $\Tor((\partial\Del_k)_+)$;
\item[({\bf E2'})] any point $\ini(\bp'_i)$ (resp., $\ini(\bp''_{i,1})$ or $\ini(\bp''_{i,2})$)
with $\bx'_i\in\ox'_2$ (resp., $\bx''_i\in\ox''_2$), dual to a
polygon $\Del_k\in P_{\cal B}(S_C)$, lies on precisely one
(possibly multiple) component of the limit curve $C_k$.
\end{enumerate}

Subtracting (\ref{espn15}) from (\ref{espn16}), we obtain
$$|\partial\Del|-2\le2r'_1+4r''_1+2r''_{2,1a}+2r''_{2,3}-2r'_{2,2}\ ,$$ whereas (\ref{espn6}) and (\ref{espn14}) lead to
$$|\partial\Del|-2\ge\sum_{\Del_k\in P(S_C)\backslash P_{\cal
B}(S_C)}\sum_{j=1}^{l_k}(B(C_{kj})-2)$$
\begin{equation}+2r'_1+4r''_1+2r'_{2,3}+2r''_{2,1a}+2r''_{2,1b}+4r''_{2,3}\
.\label{espn43}\end{equation} Comparison of the two last
inequalities results in (\ref{espn11}) together with
(\ref{espn12}), where $r'_3$ and $r''_3$ are excluded, and also
results in equalities in (\ref{espn6}), (\ref{espn8}),
(\ref{espn14}), (\ref{espn42}), and (\ref{espn15}), that is we
obtain ({\bf E3})-({\bf E5}), ({\bf E7}), and ({\bf E8}), as well
as $r''_{2,1b}=0$ (cf. (\ref{espn12})). As in the end of section
\ref{secspn3}, observe that ({\bf E5}) and ({\bf E8}) imply,
first, $r'_3=r''_3=0$, and, second, that the number of the
distinct local branches of the curves $C_k$, $1\le k\le N$,
centered along $\Tor(\sig)$, $\sig\subset\partial\Del$, is equal
to $B_+$, which in turn equals to $|\partial\Del|$, and hence all
these branches are not multiple. In particular, the above
condition ({\bf E2'}) can be replaced by ({\bf E2}) from the
preceding section, and ({\bf E1'}) reduces to
\begin{enumerate}\item[({\bf E1''})] for each polygon $\Del_k\in
P_{\cal B}(S_C)$,

- either $|(\partial\Del_k)_+|=1$, $\Del_k$ is dual to a point
from $\ox'_{2,1}\cup\ox''_{2,1a}$, the limit curve $C_k$ may have
binomial components, meeting the toric divisors $\Tor(\sig)$,
$\sig\subset(\partial\Del_k)_\perp$, and it has one real
non-multiple component, crossing the toric divisor
$\Tor((\partial\Del_k)_+)$,

- or $|(\partial\Del_k)_+|=2$, $\Del_k$ is dual to a point from
$\ox''_{2,2}$, the limit curve $C_k$ may have binomial components,
meeting the toric divisors $\Tor(\sig)$,
$\sig\subset(\partial\Del_k)_\perp$, and it has two distinct
non-multiple conjugate components, meeting the toric divisors
$\Tor((\partial\Del_k)_+)$.
\end{enumerate}

\begin{remark}\label{rspn2} In the same way as in Remark
\ref{rspn1} we decide that $\bigcup_{\Del_k\in P_{\cal
B}(S_C)}(\partial\Del_k)_+$ forms an admissible lattice path
$\gam$ satisfying conditions ($\gam${\bf1})-($\gam${\bf4}).
\end{remark}

\medskip

{\bf Case 2}. Assume that
\begin{equation}B_+=2d-d_1+n_0=|(\partial\Del)_+|+n_0,\quad n_0>0\
.\label{espn18}\end{equation}

{\it Step 1}. Similarly to Case 1, we have
$$|\partial\Del|+2n_0=2|(\partial\Del)_+|+2n_0$$ \begin{equation}\ge2(r'_{2,1}+2r'_{2,2}+r'_{2,3}+
r''_{2,1a}+2r''_{2,1b}+2r''_{2,2}+r''_{2,3}+r'_3+2r''_3)\
,\label{espn19}\end{equation} with an equality only if conditions
({\bf E1'}), ({\bf E2'}) from Case 1 hold true.

Also we refine inequality (\ref{espn6}) up to the following one.
For an edge $\sig=\Del_k\cap\Del_l$ with some $\Del_k,\Del_l\in
P(S_C)$, denote by $d(\sig)$ the absolute value of the difference
between the number of local branches of $(C_k)_\red$, centered on
$\Tor(\sig)$, and the corresponding number for $(C_l)_\red$. Then,
due to \cite{ShP}, Remark 3.4, inequality (\ref{espn6}) refines up
to
$$\sum_{\Del_k\in P_{\cal
B}(S_C)}\sum_{j=1}^{l_k}(B(C_{kj})-2)+\sum_{\Del_k\in
P(S_C)\backslash P_{\cal B}(S_C)}\sum_{j=1}^{l_k}(B(C_{kj})-2)+
\sum_{\renewcommand{\arraystretch}{0.6}
\begin{array}{c}
\scriptstyle{\sig\in E(S_C)}\\
\scriptstyle{\sig\not\subset\partial\Del}
\end{array}}d(\sig)$$
\begin{equation}\le B(\partial\Del)-2\le|\partial\Del|-2\
.\label{espn44}\end{equation}

{\it Step 2}. The key ingredient in our argument is the relation
\begin{equation}\sum_{\Del_k\in P_+(S_C)}\sum_{j=1}^{l_k}(B(C_{kj})-2)+\sum_{\renewcommand{\arraystretch}{0.6}
\begin{array}{c}
\scriptstyle{\sig\in E(S_C)}\\
\scriptstyle{\sig\subset{\cal B}_+,\ \sig\not\subset\partial\Del}
\end{array}}d(\sig)\ge n_0\ .\label{espn45}\end{equation} To
derive it, choose a generic vector $\overline v\in\R^2$, close to
$(-1,1)$, and denote by $d_{\overline v}(\Del_k,C_k)$ the
difference between the number of the local branches of
$(C_k)_\red$, centered on the divisors $\Tor(\sig)$ for the sides
$\sig$ of $\Del_k$, through which $\overline v$ enters $\Del_k$,
and the number of the local branches of $(C_k)_\red$, centered on
the divisors $\Tor(\sig)$ for the sides $\sig$ of $\Del_k$,
through which $\overline v$ leaves $\Del_k$. Clearly, for any
$\Del_k\in P(S_C)\backslash P_{\cal B}(S_C)$,
\begin{equation}\sum_{j=1}^{l_k}(B(C_{kj})-2)\ge |d_{\overline
v}(\Del_k,C_k)|\ ,\label{espn46}\end{equation} and the equality
here holds only when non-binomial components of $C_k$ are
non-multiple. Summing up inequalities (\ref{espn46}) over all
$\Del_k\in P(S_C)$, $\Del_k\subset{\cal B}_+$, and using
(\ref{espn18}), we get (\ref{espn45}).

{\it Step 3}. Subtracting (\ref{espn19}) from (\ref{espn16}), we
obtain
$$|\partial\Del|-2\le2n_0+2r'_1+4r''_1+2r''_{2,1a}+2r''_{2,3}-2r'_{2,2}\
,$$ whereas (\ref{espn14}) and (\ref{espn44}) yield
$$|\partial\Del|-2\ge B(\partial\Del)-2\ge\sum_{\Del_k\in
P(S_C)\backslash P_{\cal B}(S_C)}\sum_{j=1}^{l_k}(B(C_{kj})-2)+
\sum_{\renewcommand{\arraystretch}{0.6}
\begin{array}{c}
\scriptstyle{\sig\in E(S_C)}\\
\scriptstyle{\sig\not\subset\partial\Del}
\end{array}}d(\sig)$$ $$+2r'_1+4r''_1+2r'_{2,3}+2r''_{2,1a}+2r''_{2,1b}+4r''_{2,3}\
.$$ Together with (\ref{espn44}) and (\ref{espn45}) this implies
equalities in (\ref{espn19}), (\ref{espn44}), and (\ref{espn45})
as well as the the relations
\begin{equation}r'_{2,2}=r'_{2,3}=r''_{2,3}=r''_{2,1b}=0,\quad
B(\partial\Del)=|\partial\Del|\ ,\label{espn47}\end{equation} and
\begin{equation}\sum_{j=1}^{l_k}(B(C_{kj})-2)=d_{\overline
v}(\Del_k,C_k)=0\ \text{for all}\ \Del_k\subset{\cal B}_+\
.\label{espn50}\end{equation} Also, for any polygon
$\Del_k\subset{\cal B}_+$, we get the absence of sides,
perpendicular to ${\cal B}$ in such polygons (otherwise,
reflecting $\overline v$ with respect to the normal to ${\cal B}$,
we will break at least one relation (\ref{espn50})). In turn, the
equality conditions provide us with the properties ({\bf E1'}),
({\bf E2})-({\bf E5}), ({\bf E7}), and the following one:
\begin{enumerate}\item[({\bf E13})] for any polygon $\Del_k$ with a
side $\sig\subset{\cal B}$, the curve $C_k$ has at most two
distinct local branches centered along
$\Tor(\sig)$.\end{enumerate}

Furthermore, relation (\ref{espn18}) and the equalities in
(\ref{espn45}), (\ref{espn46}) yield that
\begin{equation}d(\sig)=0\quad\text{for all}\quad\sig\in E(S_C),\ \sig\not\subset\partial\Del\
,\label{espn48}\end{equation} in particular,
\begin{equation}\sum_{\Del_k\subset{\cal B}_+}\sum_{j=1}^{l_k}(B(C_{kj})-2)=n_0\ .\label{espn49}\end{equation}
The last equality in (\ref{espn47}), the equalities in
(\ref{espn46}), and (\ref{espn49}) tell us that
\begin{enumerate}\item[({\bf E14})] all the local branches of the curves $C_k$, $1\le
k\le N$, centered along the divisors $\Tor(\sig)$,
$\sig\subset\partial\Del$, are non-multiple, non-singular, and
transverse to $\Tor(\sig)$.\end{enumerate} Then, in particular,
the polygons $\Del_k\in P_{\cal B}(S_C)$ with
$(\partial\Del_k)_+\subset\partial\Del$ possess the property ({\bf
E1''}), introduced in Case 1.

{\it Step 4}. We shall describe the polygons $\Del_k\in
P(S_C)\backslash P_{\cal B}(S_C)$ and the respective limit curves
$C_k$.

Introduce a partial order in $P_+(S_C)$, $\Del_j\subset{\cal
B}_+$, saying that $\Del_j\prec\Del_l$ if $\Del_j\cap\Del_l=\sig$
is a common side, through which $\overline v$ goes from $\Del_j$
to $\Del_l$. Extend this partial order up to a complete one.

\begin{figure}
\begin{center}
\epsfxsize 145mm \epsfbox{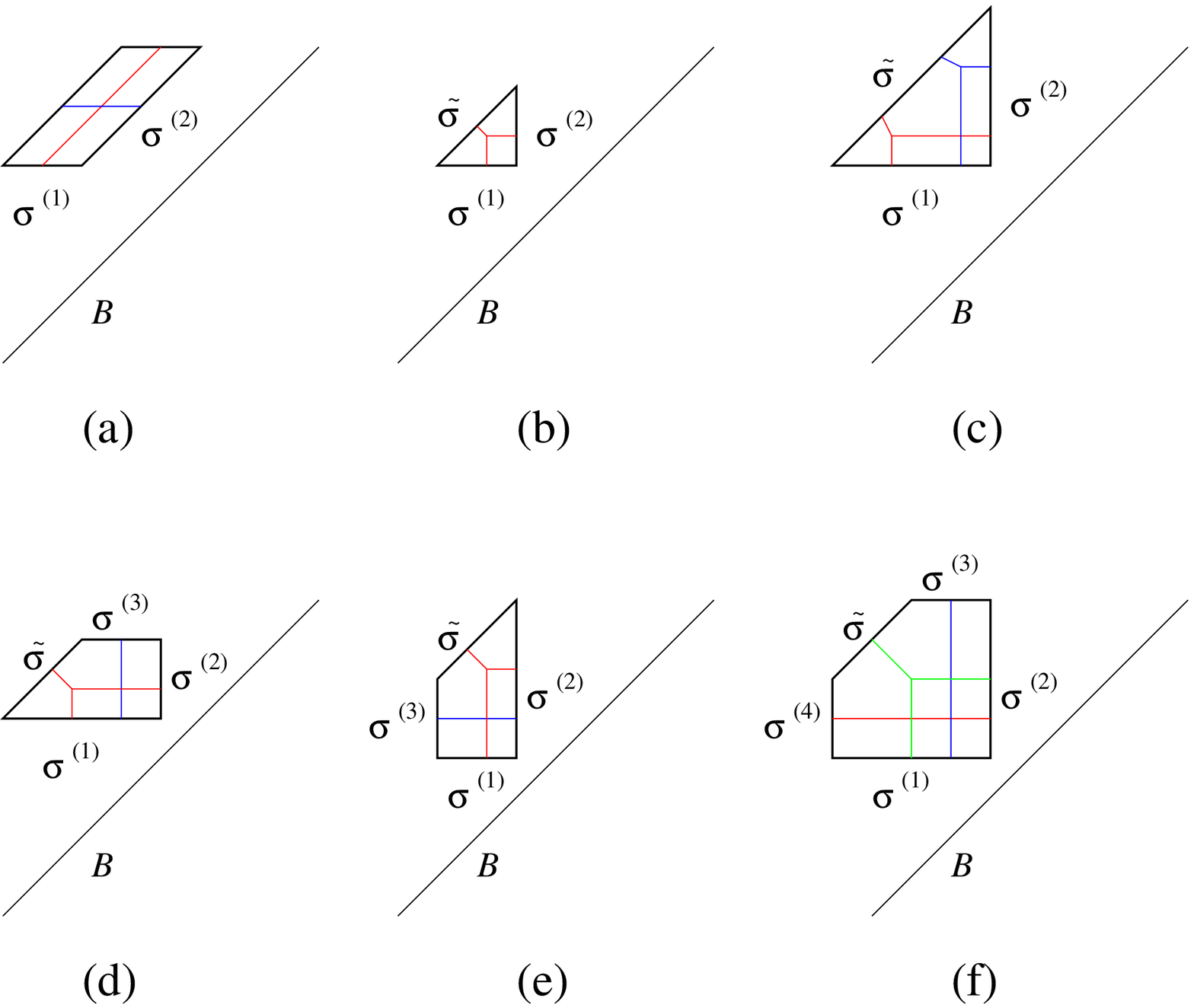}
\end{center}
\caption{Polygons in ${\cal B}_+$}\label{figsp12}
\end{figure}

Let $\Del_k\in P_+(S_C)$ be the first polygon with respect to the
order defined. Notice that $\Del_k$ has precisely two sides,
through which $\overline v$ enters $\Del_k$. Indeed, in case there
is only one such side we immediately obtain that $d_{\overline
v}(\Del_k,C_k)<0$. In case of more than two such sides, we obtain
that more than two edges of the tropical curve $A_C$, starting at
some points of $\ox$ on ${\cal B}$, merge at the vertex, dual to
$\Del_k$, which would impose a restriction to the configuration
$\ox$ in contradiction to the generality condition ($\bx${\bf1})
from section \ref{secspn6}. Let $\sig^{(1)},\sig^{(2)}$ be the
sides of $\Del_k$, through which $\overline v$ enters $\Del_k$.
Relations (\ref{espn50}), (\ref{espn48}), and properties ({\bf
E1'}), ({\bf E13}) leave the only following options for $\Del_k$
and $C_k$ (see Figure \ref{figsp12}, where colored lines designate
the components of $(C_k)_\red$ and their intersection with the
corresponding toric divisors):
\begin{itemize}\item $\Del_k$ is a parallelogram, $C_k$ splits into binomial
components - Figure \ref{figsp12}(a);
\item $\Del_k$ is a triangle,
$(C_k)_\red$ is irreducible and has only one local branch,
centered along $\Tor(\sig)$, for any side $\sig$ of $\Del_k$ -
Figure \ref{figsp12}(b);
\item $\Del_k$ is a triangle,
$(C_k)_\red$ splits into two components, these components
intersect only in $(\C^*)^2\subset\Tor(\Del_k)$, and each of them
has only one local branch, centered along $\Tor(\sig)$, for any
side $\sig$ of $\Del_k$ - Figure \ref{figsp12}(c); \item $\Del_k$
is a trapeze with a side $\sig^{(3)}$, parallel to $\sig^{(1)}$
(or $\sig^{(2)}$), the curve $(C_k)_\red$ splits into a binomial
component, crossing the divisors $\Tor(\sig^{(3)})$ and
$\Tor(\sig^{(1)})$ (resp., $\Tor(\sig^{(2)})$), and a non-binomial
components, crossing the divisors $\Tor(\sig)$,
$\sig\subset\partial\Del_k$, $\sig\ne\sig^{(3)}$, and having only
one local branch along each of these divisors - Figure
\ref{figsp12}(d) (resp., (e)); \item $\Del_k$ is a pentagon with
sides $\sig^{(3)}$, $\sig^{(4)}$, parallel to $\sig^{(1)}$,
$\sig^{(2)}$, respectively, the curve $(C_k)_\red$ splits into a
binomial component, crossing the divisors
$\Tor(\sig^{(1)}),\Tor(\sig^{(3)})$, a binomial component,
crossing the divisors $\Tor(\sig^{(2)}),\Tor(\sig^{(4)})$, and a
non-binomial component, crossing the divisors $\Tor(\sig)$,
$\sig\subset\partial\Del_k$, $\sig\ne\sig^{(3)},\sig^{(4)}$, and
having only one local branch along each of these divisors - Figure
\ref{figsp12}(f).\end{itemize}

Now take the second polygon $\Del_j\subset{\cal B}_+$ and observe
that again there are precisely two sides $\sig^{(1)},\sig^{(2)}$
of $\Del_j$, through which $\overline v$ enters $\Del_j$, and
along each of the $\Tor(\sig^{(1)}),\Tor(\sig^{(2)})$,
$(C_j)_\red$ has one or two local branches. Thus, we decide that
$\Del_j$ and $C_j$ are as shown in Figure \ref{figsp12}.
Inductively we deduce the similar conclusions for all $\Del_j\in
P_+(S_C)$. Moreover, using condition ({\bf E14}), obtained in Step
3, we can easily derive that \begin{enumerate}\item[({\bf E15})]
all $\Del_k\in P_{\cal B}(S_C)$ possess property ({\bf E1'}); all
the limit curves $C_k$, $\Del_k\in P(S_C)\backslash P_{\cal
B}(S_C)$, are reduced, they cross $\Tor(\partial\Del_k)$
transversally at their non-singular points; if $\Del_k\in
P_+(S_C)$ is as shown in Figure \ref{figsp12}(b-f), then the side
$\widetilde\sig$ is parallel to ${\cal B}$.
\end{enumerate} Taking into account condition ({\bf E5}), we
additionally obtain that $r'_3=r''_3=0$, and that $E(S_C)$ has no
edges lying on ${\cal B}$.

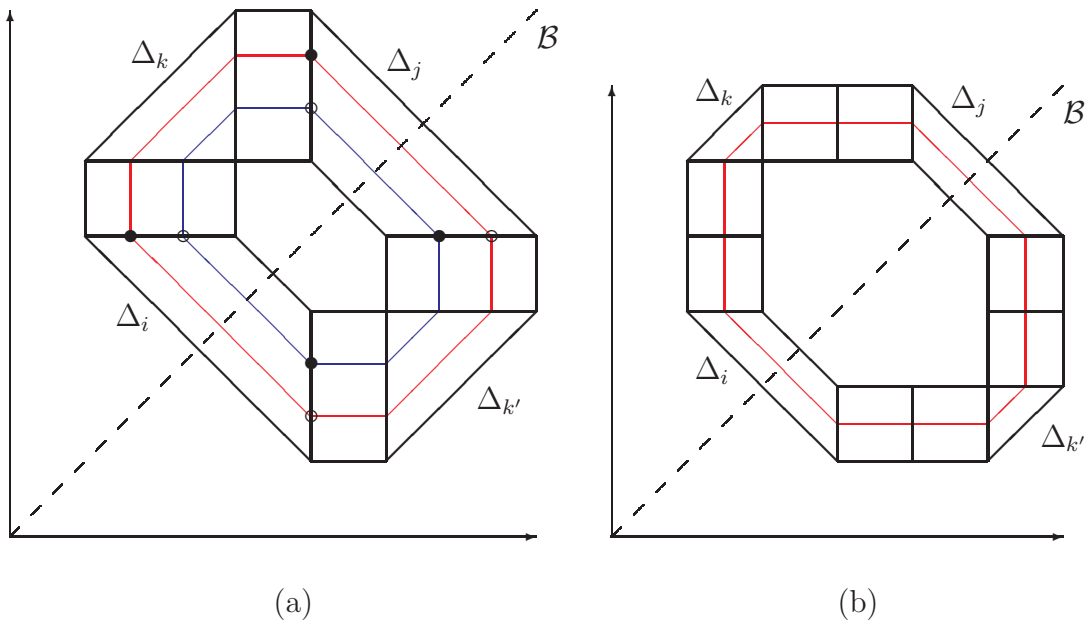
\begin{figure}
\setlength{\unitlength}{1cm}
\begin{picture}(14.5,8)(0,0)
\put(8,1){\vector(1,0){6}}\put(8,1){\vector(0,1){6}}
\put(0,1){\vector(1,0){7}}\put(0,1){\vector(0,1){7}}
{\color{Red}\put(11,2.5){\line(-1,1){1.5}}\put(11,2.5){\line(1,0){2}}
\put(13,2.5){\line(1,1){0.5}}\put(13.5,3){\line(0,1){2}}
\put(13.5,5){\line(-1,1){1.5}}\put(12,6.5){\line(-1,0){2}}
\put(10,6.5){\line(-1,-1){0.5}}\put(9.5,6){\line(0,-1){2}}
\put(1.6,5){\line(0,1){1}}\put(1.6,6){\line(1,1){1.4}}
\put(3,7.4){\line(1,0){1}}\put(4,7.4){\line(1,-1){2.4}}
\put(6.4,5){\line(0,-1){1}}\put(6.4,4){\line(-1,-1){1.4}}
\put(5,2.6){\line(-1,0){1}}\put(4,2.6){\line(-1,1){2.4}}
\color{Blue}\put(2.3,5){\line(0,1){1}}\put(2.3,6){\line(1,1){0.7}}
\put(3,6.7){\line(1,0){1}}\put(4,6.7){\line(1,-1){1.7}}
\put(5.7,5){\line(0,-1){1}}\put(5.7,4){\line(-1,-1){0.7}}
\put(5,3.3){\line(-1,0){1}}\put(4,3.3){\line(-1,1){1.7}}}
\put(1.5,4.9){$\bullet$}\put(3.9,3.2){$\bullet$}\put(3.9,7.3){$\bullet$}
\put(5.6,4.9){$\bullet$}\put(2.2,4.9){$\circ$}\put(3.9,6.6){$\circ$}
\put(3.9,2.5){$\circ$}\put(6.3,4.9){$\circ$}
\thicklines\put(9,4){\line(0,1){2}}
\put(9,6){\line(1,1){1}}\put(10,7){\line(1,0){2}}\put(12,7){\line(1,-1){2}}\put(14,5){\line(0,-1){2}}
\put(14,3){\line(-1,-1){1}}\put(13,2){\line(-1,0){2}}\put(11,2){\line(-1,1){2}}
\put(9,4){\line(1,0){1}}\put(11,3){\line(1,0){3}}\put(13,4){\line(1,0){1}}
\put(9,5){\line(1,0){1}}\put(13,5){\line(1,0){1}}\put(9,6){\line(1,0){3}}
\put(10,4){\line(0,1){3}}\put(11,2){\line(0,1){1}}\put(11,6){\line(0,1){1}}
\put(12,2){\line(0,1){1}}\put(12,6){\line(0,1){1}}\put(13,2){\line(0,1){3}}
\put(11,3){\line(-1,1){1}}\put(13,5){\line(-1,1){1}}
\put(4,2){\line(1,0){1}}\put(4,4){\line(1,0){3}}\put(5,5){\line(1,0){2}}
\put(1,5){\line(1,0){2}}\put(1,6){\line(1,0){3}}\put(3,8){\line(1,0){1}}
\put(1,5){\line(0,1){1}}\put(3,5){\line(0,1){3}}\put(4,2){\line(0,1){2}}
\put(4,6){\line(0,1){2}}\put(5,2){\line(0,1){3}}\put(7,4){\line(0,1){1}}
\put(4,2){\line(-1,1){3}}\put(4,4){\line(-1,1){1}}\put(5,5){\line(-1,1){1}}
\put(7,5){\line(-1,1){3}}\put(5,2){\line(1,1){2}}\put(1,6){\line(1,1){2}}
\put(11,0){\text{\rm (b)}}\put(3.5,0){\text{\rm
(a)}}\thinlines\dashline{0.2}(8,1)(14,7)\dashline{0.2}(0,1)(7,8)\put(14,6.5){${\cal
B}$}\put(7,7.5){${\cal
B}$}\put(9.1,6.8){$\Del_k$}\put(1.6,7.3){$\Del_k$}\put(9.1,3.1){$\Del_i$}
\put(1.4,3.8){$\Del_i$}\put(12.5,6.7){$\Del_j$}\put(5,7.2){$\Del_j$}
\put(13.7,2.2){$\Del_{k'}$}\put(6.2,2.7){$\Del_{k'}$}
\end{picture}
\caption{Forbidden subdivisions}\label{figsp13}
\end{figure}

{\it Step 5}. We claim that the triangles as shown in Figure
\ref{figsp12}(c) do not occur in $S$. Indeed, otherwise, a
triangle $\Del_k$ of this type and its symmetric with respect to
${\cal B}$ image $\Del_{k'}$ would be joined by sequences of
parallelograms with two trapezes $\Del_i,\Del_j\in P_{\cal
B}(S_C)$, having non-parallel sides of length $2$ and such that
each of the curves $C_i,C_j$ has two conjugate components,
crossing $\Tor((\partial\Del_i)_+),\Tor((\partial\Del_j)_+)$,
respectively, but then the components of the curves $C_k,C_{k'}$
and $C_i,C_j$ (together with possible binomial components of the
limit curves corresponding to parallelograms) would glue up in the
deformation $C^{(t)}$, $t\in(\C,0)$, so that the curves $C^{(t)}$,
$t\ne 0$, would have at least two handles in contradiction to
their rationality (see Figure \ref{figsp13}(a), where the gluing
components of the limit curves are designated by red and blue
lines, the bullets and circles designate pairs of conjugate points
of the limit curves on the toric divisors of $\Tor(\Del_i)$ and
$\Tor(\Del_j)$).

Denoting by $n_{0,b},n_{0,d},n_{0,e},n_{0,f}$ the number of the
polygons in $P_+(S_C)$ of the types, shown in Figure
\ref{figsp12}(b,d-f), respectively, we derive from (\ref{espn18})
and ({\bf E15}) that
\begin{equation}n_{0,b}+n_{0,d}+n_{0,e}+n_{0,f}=n_0\
.\label{espn51}\end{equation}

We also observe that the lattice path $\gam$, defined as
$\bigcup_{\Del_k\in P_{\cal B}(S_C)}(\partial\Del_k)_+$ with the
common vertices of the polygons $\Del_k\in P_{\cal B}(S_C)$ as the
points $v_i$, $1\le i<\widetilde m$, satisfies the conditions
($\gam${\bf 1})-($\gam${\bf 5}) and ($\gam${\bf 1'})-($\gam${\bf
6'}) of sections \ref{secspn1} and \ref{secspn4}. For example,
condition ($\gam${\bf 6'}) holds, since, otherwise, one would have
a polygon $\Del_k\in P_{\cal B}(S_C)$ with $(\partial\Del_k)_+$
consisting of of a segment with the projection $1/2$ to ${\cal B}$
and of a unit segment parallel to ${\cal B}$, but such a polygon
cannot split into the Minkowski sum of two polygons symmetric with
respect to ${\cal B}$, contrary to the fact that the corresponding
limit curve $C_k$ contains two conjugate components, crossing
$\Tor((\partial\Del_k)_+)$.

\begin{remark}\label{rsp1}
Notice that the condition $r''=0$ (i.e., if all the fixed points
are real) fits to the Case 1, defined by (\ref{espn42}). Indeed,
assuming in contrary (\ref{espn18}) and taking into account the
above conclusions (\ref{espn47}), ({\bf E1'}), $r'_3=0$, and
(\ref{espn51}), we derive that $(\partial\Del_i)_+$ is a unit
segment for any $\Del_i\in P_{\cal B}(S_C)$, and thus, $P_+(S_C)$
contains $n_0$ triangles like shown in Figure \ref{figsp12}(b).
Then such a triangle $\Del_k$ and its symmetric copy $\Del_{k'}$
are joined with two polygons $\Del_i,\Del_j\in P_{\cal B}(S_C)$ by
sequences of parallelograms (see, for instance, Figure
\ref{figsp13}(b)), and hence the corresponding limit curves
(designated by red lines in Figure \ref{figsp13}(b)) glue up in
the deformation $C^{(t)}$, $t\in(\C,0)$, into a non-rational curve
$C^{(t)}$, $t\ne 0$, contradicting the initial assumptions.
\end{remark}

\section{Proof of Theorems \ref{tsp1}, \ref{tsp3}, and \ref{tsp4}}

\subsection{Encoding the tropical limits}
Let us be given a surface $\Sig=\bbS^2$, $\bbS^2_{1,0}$,
$\bbS^2_{2,0}$, or $\bbS^2_{0,2}$, a respective lattice polygon
$\Del$ as shown in Figure \ref{fn2}(a-d), and a pair of
non-negative integers $r'r''$, satisfying (\ref{espn1}). Choose a
generic configuration of points $\op=\op'\cup\op''\subset\Sig(\K)$
satisfying (\ref{espn21}) and such that its valuation projection
$\ox=\ox'\cup\ox''\subset{\cal B}$ satisfies condition
($\bx${\bf1}) from section \ref{secspn6}.

For real rational curves $C\in|{\cal L}_\Del|$ on the surface
$\Sig(\K)$, passing through the configuration $\op$, we have
described possible tropical limits $(A_C,\{C_1,...,C_N\})$. We
encode them by means of the objects counted in Theorems
\ref{tsp1}, \ref{tsp3}, and \ref{tsp4}, and it is immediate from
the results of section \ref{secspn7} that they satisfy all the
conditions specified in section \ref{secspn9}. Namely,
\begin{itemize}\item the splitting ${\cal R}$ of
$r',r''$ as in (\ref{espn3}) is defined by taking $r'_2$ equal to
the number of the points $\bx'_i\in\ox'$ among the vertices of
$A_C$, and taking $r''_{2,1}$ (resp., $r''_{2,2}$) equal to the
numbers of the points $\bx''_i\in\ox''$ among the vertices of
$A_C$ such that, for the polygon $\Del_k$ in the dual subdivision
$S_C$ of $\Del$, it holds $|(\partial\Del_k)_+|=1$ (resp.,
$|(\partial\Del_k)_+|=2$); \item the broken line
$\bigcup_{\Del_k\in P_{\cal B}(S_C)}(\partial\Del_k)_+$ naturally
defines an admissible path $\gam$ with the integral points on it
$v_i$, $i=0,...,\widetilde m$, to be the endpoints of the
fragments $(\partial\Del_k)_+$, $\Del_k\in P_{\cal B}(S_C)$; \item
the intersection points of the curves $C_k$ such that $\Del_k\in
P_{\cal B}(S_C)$ with the divisors $\Tor((\partial\Del_k)_\perp)$,
form the set of the vertices of the graph $G'$, whereas the
binomial components of these curves $C_k$, crossing
$\Tor((\partial\Del_k)_\perp)$, serve as the arcs of $G'$ with
their multiplicities as weights $w(G'_i)$; in turn, picking the
vertices of $G'$, which are $\ini(\bp'_i)$ or
$\ini(\bp''_{i,1}),\ini(\bp_{i,2})$, we obtain the marking
$\overline s=\overline s'\cup\overline s''$;
\item then we take the components of the curves $C_k$,
$\Del_k\in P_{\cal B}(S)_C)$, which are not binomial crossing
$\Tor((\partial\Del_k)_\perp$, as additional vertices of the graph
$G$, joining them by arcs with the vertices of $G'$, which belong
to these components; \item finally, we take the the non-binomial
components of the curves $C_j$, $\Del_j\in P(S_C)\backslash
P_{\cal B}(S_C)$, and join them with the previously defined
vertices of $G$ as far as the corresponding components of the
limit curves either intersect, or can be connected by a sequence
of binomial components.
\end{itemize}

\begin{remark} Up to some details the graph $G$ can be viewed as a
rational parameterization of the tropical curve $A_C$.
\end{remark}

\subsection{Restoring the tropical limits}\label{secsp2} Let
a surface $\Sig$, a polygon $\Del$, a pair of nonnegative
integers $r',r''$, and a configuration $\op$ as in the preceding
section, and let ${\cal R},\gam,S,G,\overline s$ be suitable
objects from Theorems \ref{tsp1}, \ref{tsp3}, or \ref{tsp4}. We
shall describe how to recover the tropical limits of the real
rational curves on $\Sig(\K)$, compatible with the given data.

The subdivision $S$ determines the combinatorial type of the
tropical curve $A$, and to restore $A$ completely we should pick
$r'_2$ points from $\ox'$ and $r''_{2,1}+r''_{2,2}$ points from
$\ox''$ and appoint them as vertices of $A$ of the line ${\cal
B}$. This choice of $r'_2+r''_{2,1}+r''_{2,2}=\widetilde m$ points
of $\ox$ is uniquely determined by the marking $\overline s$:
namely, any $\widetilde m$ points of $\ox$ divide ${\cal B}$ into
$\widetilde m+1$ naturally ordered intervals, and the distribution
of the remaining points of $\ox'$ and $\ox''$ in these intervals
must coincide with the distribution of the values of $\overline
s'$ and $\overline s''$, respectively, in the intervals
$$\left(-\infty,\frac{1}{2}\right),\ \left(\frac{1}{2},\frac{3}{2}\right),\ ...\
\left(\widetilde m-\frac{3}{2},\widetilde m-\frac{1}{2}\right),\
\left(\widetilde m-\frac{1}{2},\infty\right)\ .$$

The tropical curve $A$ and the subdivision $S$ determine a convex
piece-wise linear function $\nu:\Del\to\R$ uniquely up to a
constant summand.

Now we pass to restoring limit curves $C_1,...,C_N$.

First, notice that, knowing the limit curves $C_k$ for all
$\Del_k\in P_{\cal B}(S)$, we uniquely define the remaining limit
curves, provided that they consist of only binomial components. In
case $\Sig=\bbS^2_{0,2}$, $r''>0$, and $n_0>0$, the limit curves
$C_j$ with $\Del_j\subset{\cal B}_+$ contain in total $n_0$
non-binomial components, which can be restored in
$\prod_{k=1}^{n_0}\Area(\widetilde\Del_k)$ ways (see the notation
of section \ref{secspn4}) due to \cite{ShP}, Lemma 3.5. For
$\Del_j\subset{\cal B}_-$ we respectively take conjugate limit
curves.

Any point $\bx'_i$ (resp., $\bx''_i$), lying inside an edge of $A$
on ${\cal B}$, defines a real point $\ini(\bp'_i)$ (resp., a pair
of conjugate points $\ini(\bp''_{i,1}),\ini(\bp''_{i,2})$) in
$\Tor(\sig)$ for the dual edge $\sig\in E(S)$. Using the weighted
graph $G'$, we can restore the respective binomial components of
the curves $C_k$, $\Del_k\in P_{\cal B}(S)$, and this can be done
in $\prod_{i=0}^{\widetilde m}
(n'_i!n''_i!\alp_i^{-1}\bet_i^{-1})$ ways (in the notation of
section \ref{secspn7}). Furthermore, for a curve $C_k$, where
$\Del_k\in P_{\cal B}(S)$ is dual to a point $\bx''_i$, its
binomial component, which crosses $\Tor((\partial\Del_k)_+)$, must
go through $\ini(\bp''_{i,1}$ or $\ini(\bp''_{i,2})$, and hence is
defined uniquely.

To restore the remaining components of $C_k$, $\Del_k\in P_{\cal
B}(S)$, and count how many solutions are there, we use

\begin{lemma}\label{lsp1} Let a lattice polygon $\Del_0$ be a triangle with a side $\sig'\perp{\cal B}$
and an opposite vertex $\sig''$, or a trapeze with sides
$\sig',\sig''\perp{\cal B}$. Let us be given generic distinct
points $z'_i\in\Tor(\sig')$, $1\le i\le p$ ($p\ge 0$), and, in
case $|\sig''|\ne 0$, generic distinct points
$z''_j\in\Tor(\sig'')$, $1\le j\le q$ ($q\ge 0$), and let $m'_i$,
$1\le i\le p$, and $m''_j$, $1\le j\le q$, be positive integers.

(1) Assume that \begin{equation}m'_1+...+m'_p=|\sig'|,\quad
m''_1+...+m''_q=|\sig''|\ ,\label{espn57}\end{equation}
$|\pr_{\cal B}(\Del_0)|=1/2$, and
$z_0\in(\C^*)^2\subset\Tor(\Del_0)$ is a generic point. Then there
is a unique rational curve $C_0\in|{\cal L}(\Del_0)|$ on the
surface $\Tor(\Del_0)$, passing through $z_0$ and satisfying
\begin{equation}C_0\cap\Tor(\sig')=\sum_{i=1}^pm'_iz'_i,\quad C_0\cap\Tor(\sig'')=
\sum_{j=1}^qm''_jz''_j\ .\label{espn56}\end{equation} Furthermore,
it is non-singular.

(2) Let $\Del_0$ be symmetric with respect to ${\cal B}$, the real
structure on the surface $\Tor(\Del_0)$ be defined by
$\conj(x,y)=(\overline y,\overline x)$, $(x,y)\in(\C^*)^2$, and
the divisors $\sum_im'_iz'_i\subset\Tor(\sig')$,
$\sum_jm''_jz''_j\subset\Tor(\sig'')$ be conjugation invariant.
\begin{enumerate}\item[(2i)] Assume that $\Del_0$ is a rectangle with a unit length
side parallel to ${\cal B}$, relation (\ref{espn57}) holds true,
$z_0\in(\R^*)^2\subset\Tor(\Del_0)$ is a generic real point, and
$C_0$ is unibranch at each point of
$C_0\cap(\Tor(\sig')\cup\Tor(\sig''))$. Then there are $2^{p+q}$
distinct real rational curves $C_0\in|{\cal L}(\Del_0)|$ on the
surface $\Tor(\Del_0)$, passing through $z_0$ and satisfying
(\ref{espn56}). Moreover, they all are non-singular along
$\Tor(\sig')$ and $\Tor(\sig'')$.
\item[(2ii)] Assume that \begin{equation}m'_1+...+m'_p=|\sig'|-(2l+1),\ l\ge 0,\quad
m''_1+...+m''_q=|\sig''|\ ,\label{espn59}\end{equation} and
$z_1,z_2\in(\C^*)^2\subset\Tor(\Del_0)$ are distinct conjugate
generic points. Then the number of real rational curves
$C_0\in|{\cal L}(\Del_0)|$ on the surface $\Tor(\Del_0)$, passing
through $z_1,z_2$, satisfying
\begin{equation}C_0\cap\Tor(\sig')=\sum_{i=1}^pm'_iz'_i+(2l+1)z'_{p+1},\quad
C_0\cap\Tor(\sig'')
=\sum_{j=1}^qm''_jz''_j\label{espn58}\end{equation} with some
$z'_{p+1}\in\Tor(\sig')$, and unibranch at each point of
$C_0\cap(\Tor(\sig')\cup\Tor(\sig''))$, is equal to $2l+1$ if
$|\pr_{\cal B}(\Del_0)|=1/2$, and is equal to $2^{p+q}$ if
$\Del_0$ is a rectangle with a unit length side parallel to ${\cal
B}$. In all the cases, $z'_{p+1}$ is a real point different from
$z'_i$, $1\le i\le p$.
\end{enumerate} Furthermore, the real rational curves in (2i), (2ii) have no solitary real nodes.
\end{lemma}

{\bf Proof}. Assume that $\pr_{\cal B}(\Del_0)=1/2$. Then a
generic curve in $|{\cal L}(\Del_0|$ is non-singular and rational.
In (1) and (2i) we impose a complete set of generic linear
conditions to $C_0$, which gives the required uniqueness of $C_0$.
In case (2ii), we can define $C_0$ by an equation
$$F(x,y):=\sqrt{-1}(x\overline\lam-\lam y)^{2l+1}P(x,y)+Q(x,y)=0\ ,$$ where $P,Q$ are
given generic homogeneous polynomials, satisfying
$$\overline{P(x,y)}=P(\overline y,\overline x),\ \overline{Q(x,y)}=Q(\overline y,\overline x),
\quad \deg Q-\deg P=2l,$$ and $(\lam,\overline\lam)\in\C
P^1\simeq\Tor(\sig')$ are unknown coordinates of $z'_{p+1}$. which
can be found from the relation $F(z_1)=0$ that gives us $2l+1$
solutions.

Assume that $\Del_0$ is a rectangle with a unit length side
parallel to ${\cal B}$.

Under conditions (\ref{espn57}), (\ref{espn56}), a curve
$C_0\backslash\{z_0\}$ possesses a parameterization $\Pi:\C \to
C_0$ such that
$$\Pi(\lam'_i)=z'_i,\ 1\le i\le p,\quad\Pi(\lam''_j)=z''_j,\
1\le j\le q,\quad\Pi(\tau)\in\Tor((\Del_0)_+)\ ,$$ for some
conjugation invariant tuples
$(\lam'_1,...,\lam'_p),(\lam''_1,...,\lam''_q)\subset\C$ and
$\tau\in\C\backslash\R$. We can also assume that a pair of real or
conjugate values among $\lam'_1,...,\lam''_q$ are fixed. So, $\Pi$
is given by
$$x(t)=\xi\cdot\frac{t-\tau}{t-\overline\tau}\cdot
\prod_{i=1}^p(t-\lam'_i)^{m'_i}\cdot\left(\prod_{j=1}^q(t-\lam''_j)^{m''_j}\right)^{-1}\
,$$
$$y(t)=\overline\xi\cdot\frac{t-\overline\tau}{t-\tau}\cdot
\prod_{i=1}^p(t-\lam'_i)^{m'_i}\cdot\left(\prod_{j=1}^q(t-\lam''_j)^{m''_j}\right)^{-1}\
,$$ where $z_0=(\xi,\overline\xi)\in(\C^*)^2$. Evaluating these
equations at the two fixed parameters among
$\lam'_1,...,\lam''_q$, we obtain two equations like
$$\left(\frac{a'-\tau}{b'-\overline\tau}\right)^2=c',\quad\left(\frac{a''-\tau}{b''-\overline\tau}\right)^2=c''$$
with either two real generic triples $(a',b',c')$,
$(a'',b'',c'')$, or two generic complex conjugate triples
$(a',b',c')$, $(a'',b'',c'')$, and in both the cases we obtain
four values of $\tau$. Next, for any other parameter $\lam$ among
$\lam'_1,...,\lam''_q$, we obtain an equation like
$$\left(\frac{\lam-\tau}{\lam-\overline\tau}\right)^2=a$$ with some generic (real or complex) $a$, which has two
solutions. That is we are done in case (2i).

Under conditions (\ref{espn59}), (\ref{espn58}), there is a real
point among $z'_1,...,z'_p,z''_1,...,z''_q$. Without loss of
generality we can assume that this is $z''_1$. A curve
$C_0\backslash\{z''_1\}$ possesses a parameterization $\Pi:\C\to
C_0$ such that $$\Pi(\lam'_i)=z'_i,\ 1\le i\le
p,\quad\Pi(\lam''_j)=z''_j,\ 2\le j\le q,\quad\Pi(\sqrt{-1})=z_1\
,$$ with unknown parameters
$\lam'_1,...,\lam'_p,\lam''_2,...,\lam''_q$. We then can write
down $\Pi$ as
$$x(t)=\xi\cdot\frac{t-\tau}{t-\overline\tau}\cdot
\prod_{i=1}^p(t-\lam'_i)^{m'_i}\cdot(t-\lam_0)^{2l+1}\cdot\left(\prod_{j=2}^q(t-\lam''_j)^{m''_j}\right)^{-1}\
,$$
$$y(t)=\overline\xi\cdot\frac{t-\overline\tau}{t-\tau}\cdot
\prod_{i=1}^p(t-\lam'_i)^{m'_i}\cdot(t-\lam_0)^{2l+1}\cdot\left(\prod_{j=2}^q(t-\lam''_j)^{m''_j}\right)^{-1}\
,$$ where $z''_1=(\xi,\overline\xi)\in\C P^1\simeq\Tor(\sig'')$
and $\Pi(\lam_0)=z'_{p+1}$. From the relation $\Pi(\sqrt{-1})=z_1$
we derive an equation like
$$\frac{x(\sqrt{-1})}{y(\sqrt{-1})}=\frac{\xi}{\overline\xi}
\cdot\left(\frac{\sqrt{-1}-\tau}{\sqrt{-1}-\overline\tau}\right)^2=a$$
with a given generic $a\in\C^*$ and obtain two values of $\tau$.
For each $\lam'_i$, $1\le i\le p$ (and similarly, for each
$\lam''_j$, $2\le j\le q$), we obtain an equations in the form
$$\frac{x(\lam'_i)}{y(\lam'_i)}=\frac{\xi}{\overline\xi}\cdot\left(\frac{\lam'_i
-\tau}{\lam'_i-\overline\tau}\right)^2=a'_i$$ with $a'_i\ne 0$,
which gives two values for $\lam'_i$. Finally, we extract the
$x$-coordinate relation from $\Pi(\sqrt{-1})=z_1$\ :
$$\xi\cdot\frac{\sqrt{-1}-\tau}{\sqrt{-1}-\overline\tau}\cdot
\prod_{i=1}^p(\sqrt{-1}-\lam'_i)^{m'_i}\cdot(\sqrt{-1}
-\lam_0)^{2l+1}\cdot\left(\prod_{j=2}^q(\sqrt{-1}-\lam''_j)^{m''_j}\right)^{-1}=a_1\
,$$ where $a_1$ is the first coordinate of $z_1$, and obtain
$2l+1$ (real) solutions for $\lam_0$ (recall that $\xi$ is defined
up to a nonzero real factor). Claim (2ii) follows.

To describe the real nodes of $C_0$ in the case of the rectangular
$\Del_0$, we consider the double cover $C_0\to\C P^1$ defined by
the pencil $\alp x+\bet y=0$, $(\alp,\bet)\in\C P^1$, with two
imaginary conjugate ramification points $(1,0)$ and $(0,1)$. The
above parameterizations show that $\dim C_0(\R)=1$, and hence
$C_0$ has no solitary real nodes, since otherwise there would be
real ramification points. \proofend

With this statement we finally obtain that the number of
admissible tropical limits, corresponding to all odd weights
$w(G'_j)$, $j>2r''_1$, in the combinatorial data, is given by the
right hand side of (\ref{espn60}) if $r''=0$, or the absolute
value of the right hand side of (\ref{espn61}), (\ref{espn62})
with the reduced the term $\prod_{j=1}^{r''_1}(w(G'_{2j}))^2$, if
$r''>0$.

\subsection{Computation of Welschinger invariants} To recover real
rational curves $C\in|{\cal L}_\Del|$ on the surface $\Sig(\K)$,
passing through $\op$, we make use of the patchworking theory from
\cite{ShP}, section 5, and \cite{ShW}, section 3.

Namely, a tropical limit $(A,\{C_1,...,C_N\})$ as constructed
above, can be completed by deformation patterns, which are
associated with the components of the graph $G'$ having weight
$>1$ and are represented by rational curves with Newton triangles
$\conv\{(0,0),(0,2),(m,1)\}$, $m$ being the weight of the
corresponding component of $G'$ (see \cite{ShP}, sections 3.5,
3.6). To obtain a real rational curve $C\in|{\cal L}_\Del|$, we
have to choose \begin{itemize}\item a pair of conjugate
deformation patterns for each pair of components of $G'$,
corresponding to two conjugate imaginary points of the set $\Phi$
(the set of intersection points of the limit curves with the
divisors $\Tor(\sig)$, $\sig\in E(S)$, $\sig\perp{\cal B}$, as
defined in section \ref{secspn3}); \item a real deformation
pattern for each component of $G'$ with weight $>1$, corresponding
to a real point of $\Phi$.\end{itemize} By \cite{ShP}, Lemma 3.9,
\begin{itemize}\item for any component of $G'$ with weight $m>1$,
there are $m$ (complex) deformation patterns; \item for a
component of $G'$ with even weight, corresponding to a real point
of $\Phi$, there are no real deformation patterns, or there are
two real deformation patterns, one having an odd number of real
solitary nodes, and the other having no real solitary nodes, \item
for a component of $G'$ with odd weight, corresponding to a real
point of $\Phi$, there is precisely one real deformation pattern
and it has an even number of real solitary nodes.\end{itemize} For
each choice of an admissible tropical limit and suitable set of
deformation patterns, Theorem 5 from \cite{ShP} produces a family
of real rational curves $C\in|{\cal L}_\Del|$ on $\Sig(\K)$, which
smoothly depends on $r'+2r''=|\partial\Del|-1$ parameters.
Moreover, any curve in the family has the same number of real
solitary nodes, equal to the total number of real solitary nodes
of the real limit curves and real deformation patterns.

Then we fix the parameters by imposing the condition to pass
through $\op$. By \cite{ShW}, section 3.2, a point $\bp\in\op$
such that $\bx=\val(\bp)$ is a vertex of $A$ and
$\ini(\bp)\in(\C^*)^2\subset\Tor(\Del_k)$ for some $\Del_k\in
P_{\cal B}(S)$, gives one smooth relation to the parameters (see
\cite{ShW}, formula (3.7)). By \cite{ShP}, section 5.4, a point
$\bp\in\op$ such that $\bx=\val(\bp)$ lies inside an edge of $A$
and $\ini(\bp)\in\Tor(\sig)$ for some $\sig\in E(S)$,
$\sig\subset\Del_k$, gives a choice of $m$ smooth equations on the
parameters, where $m$ is the intersection number of $\Tor(\sig)$
with $C_k$ at $\ini(\bp)$ (see \cite{ShP}, formula (5.4.26)). In
the latter case, assume that $\ini(\bp)$ is real. Then among $m$
relations to parameters there is one real, if $m$ is odd, and
there are zero or two real relations if $m$ is even.

Choosing one relation for each point of $\op$, we obtain a
transverse system (see \cite{ShW}, section 3.2), and hence one
real rational curve $C\in|{\cal L}_\Del|$ passing through $\op$.

We notice that whenever an even weight $m$ is assigned to a
component of $G'$, corresponding to a real point of the set
$\Phi$, then the real rational curves in $|{\cal L}_\Del|$ with a
given admissible tropical limit, either do not exist, or can be
arranged in pairs with opposite Welschinger numbers due to the
choice of two deformation patterns with distinct parity of the
number of real solitary nodes (cf. \cite{ShP}, proof of
Proposition 6.1), and hence all these real rational curves do not
contribute to the Welschinger invariant. If all the weights of the
components of $G'$, corresponding to real points in $\Phi$, are
odd, then each admissible tropical limit bears
$\prod_{j=1}^{r''_1}(w(G'_{2j}))^2$ real rational curves in
$|{\cal L}_\Del|$ due to the choice of
$\prod_{j=1}^{r''_1}w(G'_{2j})$ suitable conjugation invariant
collections of deformation patterns and the choice of
$\prod_{j=1}^{r''_1}w(G'_{2j})$ conjugation invariant collections
of relations imposed by the condition $C\supset\op$. All these
curves have the same number of real solitary nodes, and the parity
of these numbers coincides with the parity of the total number of
the intersection points of the imaginary conjugate components of
the limit curves $C_k$ with $\Del_k\in P_{\cal B}(S)$. So, the
Welschinger numbers of the real rational curves in count are equal
to $(-1)^a$ as defined in formulas (\ref{espn61}), (\ref{espn62}).

{\it Address}: School of Mathematical Sciences, Tel Aviv
University, Ramat Aviv, 69978 Tel Aviv, Israel.

{\it E-mail}: shustin@post.tau.ac.il

\end{document}